\documentclass[ejs,preprint]{imsart}

\arxiv{0901.1342}

\usepackage{amsmath,amssymb,amsfonts,amscd,graphics}
\usepackage[square]{natbib}
\usepackage{hyperref}

\newcommand{\Expect}[1]{{\mathbf{E}}\left[#1\right]}
\DeclareMathOperator*{\essinf}{ess\ inf}
\DeclareMathOperator{\supp}{supp}
\newcommand{\myexp}[1]{\ensuremath{\exp}{\left\{ #1\right\}}}
\newcommand{\Posterior}{\ensuremath{\Pi}}
\newcommand{\posteriordensity}{\ensuremath{\pi}}
\newcommand{\ParameterSpace}{\ensuremath{\Theta}}
\newcommand{\ParameterSigmaField}{\ensuremath{\mathcal{T}}}
\newcommand{\parametervalue}{\ensuremath{\theta}}
\newcommand{\popavgt}[2]{\Posterior_{#2} \left( #1\right)}
\newcommand{\prioravg}[1]{\popavgt{ #1}{0}}
\newcommand{\Filtration}[1]{\ensuremath{\sigma\left(X_1^{#1}\right)}}
\newcommand{\ObservableSpace}{\ensuremath{\Xi}}
\newcommand{\ObservableSigmaField}{\ensuremath{\mathcal{X}}}
\newcommand{\LikelihoodRatio}{\ensuremath{R}}
\newcommand{\CondLike}{\ensuremath{L}}
\newcommand{\InfinitySet}{\ensuremath{I}}
\newcommand{\truedensity}{p}
\newcommand{\TrueMeasure}{P}
\newcommand{\ModelMeasure}{F_{\parametervalue}}
\newcommand{\modeldensity}{f_{\parametervalue}}
\newcommand{\altmodeldensity}{f_{\parametervalue^{\prime}}}
\newcommand{\PredictiveMeasure}{F_{\Posterior}}
\newcommand{\divrate}{h}
\newcommand{\Action}{J}
\newcommand{\convtime}{\tau}
\newtheorem{assumption}{Assumption}
\newtheorem{theorem}{Theorem}
\newtheorem{lemma}{Lemma}
\newtheorem{corollary}{Corollary}

\begin{document}
\begin{frontmatter}

\title{Dynamics of Bayesian Updating with Dependent Data and Misspecified Models}
\runtitle{Dynamics of Bayesian Updating}
\begin{aug}
\author{Cosma Rohilla Shalizi\ead[label=e1]{cshalizi@cmu.edu}}
%\affiliation{Statistics Department, Carnegie Mellon University}
\address{Statistics Department\\ Carnegie Mellon University\\ 5000 Forbes Avenue\\ Pittsburgh, PA\\ 15213-3890 USA \printead{e1}}
\runauthor{C. R. Shalizi}
\end{aug}

\begin{abstract}
  Much is now known about the consistency of Bayesian updating on
  infinite-dimensional parameter spaces with independent or Markovian data.
  Necessary conditions for consistency include the prior putting enough weight
  on the correct neighborhoods of the data-generating distribution; various
  sufficient conditions further restrict the prior in ways analogous to
  capacity control in frequentist nonparametrics. The asymptotics of Bayesian
  updating with mis-specified models or priors, or non-Markovian data, are far
  less well explored.  Here I establish sufficient conditions for posterior
  convergence when all hypotheses are wrong, and the data have complex
  dependencies.  The main dynamical assumption is the asymptotic equipartition
  (Shannon-McMillan-Breiman) property of information theory.  This, along with
  Egorov's Theorem on uniform convergence, lets me build a sieve-like structure
  for the prior.  The main statistical assumption, also a form of capacity
  control, concerns the compatibility of the prior and the data-generating
  process, controlling the fluctuations in the log-likelihood when averaged over
  the sieve-like sets.  In addition to posterior convergence, I derive a kind
  of large deviations principle for the posterior measure, extending in some
  cases to rates of convergence, and discuss the advantages of predicting using
  a combination of models known to be wrong.  An appendix sketches connections
  between these results and the replicator dynamics of evolutionary theory.
\end{abstract}



\end{frontmatter}

\section{Introduction}

The problem of the convergence and frequentist consistency of Bayesian learning
goes as follows.  We encounter observations $X_1, X_2, \ldots$, which we would
like to predict by means of a family $\ParameterSpace$ of models or hypotheses
(indexed by $\parametervalue$).  We begin with a prior probability distribution
$\Posterior_0$ over $\ParameterSpace$, and update this using Bayes's rule, so
that our distribution after seeing $X_1, X_2, \ldots X_t \equiv X_1^t$ is
$\Posterior_t$.  If the observations come from a stochastic process with
infinite-dimensional distribution $\TrueMeasure$, when does $\Posterior_t$
converge $\TrueMeasure$-almost surely?  What is the rate of convergence?  Under
what conditions will Bayesian learning be consistent, so that $\Posterior_t$
doesn't just converge but its limit is $\TrueMeasure$?

Since the Bayesian estimate is really the whole posterior probability
distribution $\Posterior_t$ rather than a point or set in $\ParameterSpace$,
consistency becomes concentration of $\Posterior_t$ around $\TrueMeasure$.  One
defines some sufficiently strong set of neighborhoods of $\TrueMeasure$ in the
space of probability distributions on $X_1^{\infty}$, and says that
$\Posterior_t$ is consistent when, for each such neighborhood $N$,
$\lim_{t\rightarrow\infty}{\Posterior_t N} = 1$.  When this holds, the
posterior increasingly approximates a degenerate distribution centered at the
truth.

The greatest importance of these problems, perhaps, is their bearing on the
objectivity and reliability of Bayesian inference; consistency proofs and
convergence rates are, as it were, frequentist licenses for Bayesian practices.
Moreover, if Bayesian learners starting from different priors converge rapidly
on the same posterior distribution, there is less reason to worry about the
subjective or arbitrary element in the choice of priors.  (Such ``merger of
opinion'' results \citep{Blackwell-Dubins-merging} are also important in
economics and game theory \citep{Chamley-herds}.)  Recent years have seen
considerable work on these problems, especially in the non-parametric setting
where the model space $\Theta$ is infinite-dimensional
\citep{Ghosh-Ramamoorthi}.

Pioneering work by Doob \citep{Doob-bayesian-consistency}, using elegant
martingale arguments, established that when any consistent estimator exists,
and $\TrueMeasure$ lies in the support of $\Posterior_0$, the set of sample
paths on which the Bayesian learner fails to converge to the truth has prior
probability zero.  (See \citep{Choi-Ramamoorthi-posterior-consistency} and
\citep{Bayesian-consistency-for-stationary-models} for extensions of this result
to non-IID settings, and also the discussion in
\citep{Schervish-theory-of-stats,Earman-on-Bayes}.)  This is not, however,
totally reassuring, since $\TrueMeasure$ generally also has prior probability
zero, and it would be unfortunate if these two measure-zero sets should happen
to coincide.  Indeed, Diaconis and Freedman established that the consistency of
Bayesian inference depends crucially on the choice of prior, and that even very
natural priors can lead to inconsistency (see \citep{Diaconis-Freedman-1986} and
references therein).

Subsequent work, following a path established by Schwartz
\citep{Schwartz-on-Bayes-procedures}, has shown that, no matter what the true
data-generating distribution $\TrueMeasure$, Bayesian updating converges along
$\TrueMeasure$-almost-all sample paths, provided that (a) $\TrueMeasure$ is
contained in $\ParameterSpace$, (b) every Kullback-Leibler neighborhood in the
$\ParameterSpace$ has some positive prior probability (the ``Kullback-Leibler
property''), and (c) certain restrictions hold on the prior, amounting to
versions of capacity control, as in the method of sieves or structural risk
minimization.  These contributions also make (d) certain dynamical assumptions
about the data-generating process, most often that it is IID
\citep{Barron-Schervish-Wasserman,Ghosal-et-al-consistency-issues,Walker-new-approaches-to-bayesian-consistency}
(in this setting, \citep{Ghosal-Ghosh-van-der-Vaart} and
\citep{Shen-Wasserman-rates} in particular consider convergence rates),
independent non-identically distributed
\citep{Choi-Ramamoorthi-posterior-consistency,Ghosal-van-der-Vaart-posterior-convergence-non-iid},
or, in some cases, Markovian
\citep{Ghosal-Tang-bayesian-consistency-for-Markov,Ghosal-van-der-Vaart-posterior-convergence-non-iid};
\citep{Choudhuri-Ghosal-Roy-bayesian-spectra} and
\citep{Roy-Ghosal-Rosenberger-sequential-bayes} work with spectral density
estimation and sequential analysis, respectively, again exploiting specific
dynamical properties.

For mis-specified models, that is settings where (a) above fails, important
early results were obtained by Berk
\citep{Berk-limiting-behavior-of-posterior,Berk-consistency} for IID data,
albeit under rather strong restrictions on likelihood functions and parameter
spaces, showing that the posterior distribution concentrates on an ``asymptotic
carrier'', consisting of the hypotheses which are the best available
approximations, in the Kullback-Leibler sense, to $\TrueMeasure$ within the
support of the prior.  More recently, \citep{Kleijn-van-der-Vaart},
\citep{Zhang-from-epsilon-to-KL} and \citep{Lian-rates-under-misspecification}
have dealt with the convergence of non-parametric Bayesian estimation for IID
data when $\TrueMeasure$ is not in the support of the prior, obtaining results
similar to Berk's in far more general settings, extending in some situations to
rates of convergence.  All of this work, however, relies on the dynamical
assumption of an IID data-source.

This paper gives sufficient conditions for the convergence of the posterior
without assuming (a), and substantially weakening (c) and (d).  Even if one
uses non-parametric models, cases where one knows that the true data generating
process is exactly represented by one of the hypotheses in the model class are
scarce.  Moreover, while IID data can be produced, with some trouble and
expense, in the laboratory or in a well-conducted survey, in many applications
the data are not just heterogeneous and dependent, but their heterogeneity and
dependence is precisely what is of interest.  This raises the question of what
Bayesian updating does when the truth is not contained in the support of the
prior, and observations have complicated dependencies.

To answer this question, I first weaken the dynamical assumptions to the
asymptotic equipartition property (Shannon-McMillan-Breiman theorem) of
information theory, i.e., for each hypothesis $\parametervalue$, the
log-likelihood per unit time converges almost surely.  This log-likelihood per
unit time is basically the growth rate of the Kullback-Leibler divergence
between $\TrueMeasure$ and $\parametervalue$, $\divrate(\parametervalue)$.  As
observations accumulate, areas of $\ParameterSpace$ where
$\divrate(\parametervalue)$ exceeds its essential infimum
$\divrate(\ParameterSpace)$ tend to lose posterior probability, which
concentrates in divergence-minimizing regions.  Some additional conditions on
the prior distribution are needed to prevent it from putting too much weight
initially on hypotheses with high divergence rates but slow convergence of the
log-likelihood.  As the latter assumptions are strengthened, more and more can
be said about the convergence of the posterior.

Using the weakest set of conditions (Assumptions
\ref{assumption:measurable-likelihood}--\ref{assumption:relative-AEP-holds}),
the long-run exponential growth rate of the posterior density at
$\parametervalue$ cannot exceed $\divrate(\Theta) - \divrate(\theta)$ (Theorem
\ref{thm:upper-bound-on-long-run-fitness}).  Adding Assumptions
\ref{assumption:best-divergence-is-finite}--\ref{assumption:sufficiently-rapid-convergence}
to provide better control over the integrated or marginal likelihood
establishes (Theorem \ref{theorem:long-run-fitness}) that the long-run growth
rate of the posterior density is in fact $\divrate(\Theta) - \divrate(\theta)$.
One more assumption (\ref{assumption:good-sets-are-good-even-for-subsets}) then
lets us conclude (Theorem \ref{thm:convergence-of-measure}) that the posterior
distribution converges, in the sense that, for any set of hypotheses $A$, the
posterior probability $\Posterior_t(A) \rightarrow 0$ unless the essential
infimum of $\divrate(\theta)$ over $A$ equals $\divrate(\Theta)$.  In fact, we
then have a kind of large deviations principle for the posterior measure
(Theorem \ref{thm:ldp}), as well as a bound on the generalization ability of
the posterior predictive distribution (Theorem
\ref{thm:generalization-performance}).  Convergence rates for the posterior
(Theorem \ref{thm:rate-of-convergence}) follow from the combination of the
large deviations result with an extra condition related to assumption
\ref{assumption:sufficiently-rapid-convergence}.  Importantly, Assumptions
\ref{assumption:best-divergence-is-finite}--\ref{assumption:good-sets-are-good-even-for-subsets},
and so the results following from them, involve {\em both} the prior
distribution and the data-generating process, and require the former to be
adapted to the latter.  Under mis-specification, it does not seem to be
possible to guarantee posterior convergence by conditions on the prior alone,
at least not with the techniques used here.

For the convenience of reader, the development uses the usual statistical
vocabulary and machinery.  In addition to the asymptotic equipartition
property, the main technical tools are on the one hand Egorov's theorem from
basic measure theory, which is used to construct a sieve-like sequence of sets
on which log-likelihood ratios converge uniformly, and on the other hand
Assumption \ref{assumption:sufficiently-rapid-convergence} bounding how long
averages over these sets can remain far from their long-run limits.  The latter
assumption is crucial, novel, and, in its present form, awkward to check; I
take up its relation to more familiar assumptions in the discussion.  It may be
of interest, however, that the results were first found via an apparently-novel
analogy between Bayesian updating and the ``replicator equation'' of
evolutionary dynamics, which is a formalization of the Darwinian idea of
natural selection.  Individual hypotheses play the role of distinct replicators
in a population, the posterior distribution being the population distribution
over replicators and fitness being proportional to likelihood.  Appendix
\ref{sec:updating-is-replication} gives details.

\section{Preliminaries and Notation}

Let $(\Omega, \mathcal{F}, P)$ be a probability space, and $X_1, X_2, \ldots$,
for short $X_1^\infty$, be a sequence of random variables, taking values in the
measurable space $(\ObservableSpace,\ObservableSigmaField)$, whose
infinite-dimensional distribution is $\TrueMeasure$.  The natural filtration of
this process is $\Filtration{t}$.  The only dynamical properties are those
required for the Shannon-McMillan-Breiman theorem (Assumption
\ref{assumption:relative-AEP-holds}); more specific assumptions such as
$\TrueMeasure$ being a product measure, Markovian, exchangeable, etc., are not
required.  Unless otherwise noted, all probabilities are taken with respect to
$\TrueMeasure$, and $\Expect{\cdot}$ always means expectation under that
distribution.

Statistical hypotheses, i.e., distributions of processes adapted to
$\Filtration{t}$, are denoted by $\ModelMeasure$, the index $\parametervalue$
taking values in the hypothesis space, a measurable space
$(\ParameterSpace,\ParameterSigmaField)$, generally infinite-dimensional.  For
convenience, assume that $\TrueMeasure$ and all the $\ModelMeasure$ are
dominated by a common reference measure, with respective densities
$\truedensity$ and $\modeldensity$.  I do not assume that
$\TrueMeasure\in\ParameterSpace$, still less that $\TrueMeasure
\in\supp{\Posterior}_0$ --- i.e., quite possibly all of the available
hypotheses are false.

We will study the evolution of a sequence of probability measures
$\Posterior_t$ on $(\ParameterSpace,\ParameterSigmaField)$, starting with a
{\em non-random} prior measure $\Posterior_0$.  (A filtration on
$\ParameterSpace$ is not needed; the measures $\Posterior_t$ change but not the
$\sigma$-field $\ParameterSigmaField$.)  Assume all $\Posterior_t$ are
absolutely continuous with respect to a common reference measure, with
densities $\posteriordensity_t$.  Expectations with respect to these measures
will be written either as explicit integrals or de Finetti style,
$\popavgt{f}{t} = \int{f(\parametervalue) d\Posterior_t(\parametervalue)}$;
when $A$ is a set, $\popavgt{f A}{t} = \popavgt{f \mathbf{1}_A}{t} =
\int_{A}{f(\parametervalue) d\Posterior_t(\parametervalue)}$.

Let $\CondLike_t(\parametervalue)$ be the conditional likelihood of $x_t$ under
$\parametervalue$, i.e., $\CondLike_t(\parametervalue) \equiv \modeldensity(X_t
= x_t|X^{t-1}_{1} = x^{t-1}_1)$, with $\CondLike_0 = 1$.  The integrated conditional
likelihood is $\popavgt{\CondLike_t}{t}$.
Bayesian updating of course means that, for any $A \in \ParameterSigmaField$,
\[
\popavgt{A}{t+1} = \frac{\popavgt{\CondLike_{t+1} A}{t}}{\popavgt{\CondLike_{t+1}}{t}}
\]
or, in terms of the density,
\[
\posteriordensity_{t+1}(\parametervalue) = \frac{\CondLike_{t+1}(\parametervalue) \posteriordensity_t(\parametervalue)}{\popavgt{\CondLike_{t+1}}{t}}
\]

It will also be convenient to express Bayesian updating in terms of the prior
and the total likelihood:
\[
\popavgt{A}{t} = \frac{\int_{A}{d\Posterior_0(\parametervalue) \modeldensity(x_1^t)}}{\int_{\ParameterSpace}{d\Posterior_0(\parametervalue) \modeldensity(x_1^t)}} = \frac{\int_{A}{d\Posterior_0(\parametervalue) \frac{\modeldensity(x_1^t)}{\truedensity(x_1^t)}}}{\int_{\ParameterSpace}{d\Posterior_0(\parametervalue) \frac{\modeldensity(x_1^t)}{\truedensity(x_1^t)}}}
 =  \frac{\prioravg{\LikelihoodRatio_t A}}{\prioravg{\LikelihoodRatio_t}}
\]
where $\LikelihoodRatio_t(\parametervalue) \equiv
\frac{\modeldensity(x_1^t)}{\truedensity(x_1^t)}$ is the ratio of model
likelihood to true likelihood.  (Note that $0 < \truedensity(x_1^t) < \infty$
for all $t$ $\TrueMeasure$-a.s.)  Similarly,
\[
\posteriordensity_t(\parametervalue) = \posteriordensity_0(\parametervalue) \frac{\LikelihoodRatio_t(\parametervalue)}{\prioravg{\LikelihoodRatio_t}}
\]

The one-step-ahead predictive distribution of the hypothesis $\parametervalue$
is given by $\ModelMeasure\left(X_t|\Filtration{t-1}\right)$, with the
convention that $t=1$ gives the marginal distribution of the first observation.
Abbreviate this by $\ModelMeasure^t$.  Similarly, let $\TrueMeasure^t \equiv
\TrueMeasure\left(X_t|\Filtration{t-1}\right)$; this is the best probabilistic
prediction we could make, did we but know $\TrueMeasure$
\citep{Knight-predictive-view}.  The posterior predictive distribution is given
by mixing the individual predictive distributions with weights given by the
posterior:
\[
\PredictiveMeasure^t \equiv \int_{\ParameterSpace}{\ModelMeasure^t d\Posterior_t(\parametervalue)}
\]

\paragraph{Remark on the topology of $\ParameterSpace$ and on $\ParameterSigmaField$}
The hope in studying posterior convergence is to show that, as $t$ grows, with
higher and higher ($\TrueMeasure$) probability, $\Posterior_t$ concentrates
more and more on sets which come closer and closer to $\TrueMeasure$.  The
tricky part here is ``closer and closer'': points in $\ParameterSpace$
represent infinite-dimensional stochastic process distributions, and the
topology of such spaces is somewhat odd, and irritatingly abrupt, at least
under the more common distances.  Any two ergodic measures are either equal or
have completely disjoint supports \citep{Gray-ergodic-properties}, so that the
Kullback-Leibler divergence between distinct ergodic processes is always
infinity (in both directions), and the total variation and Hellinger distances
are likewise maximal.  Most previous work on posterior consistency has
restricted itself to models where the infinite-dimensional process
distributions are formed by products of fixed-dimensional base distributions
(IID, Markov, etc.), and in effect transferred the usual metrics' topologies
from these finite-dimensional distributions to the processes.  It {\em is}
possible to define metrics for general stochastic processes
\citep{Gray-ergodic-properties}, and if readers like they may imagine that
$\ParameterSigmaField$ is a Borel $\sigma$-field under some such metric.  This
is not necessary for the results presented here, however.

\subsection{Example} The following example will be used to illustrate the
assumptions (\S \ref{sec:verification-of-assumptions} and Appendix
\ref{appendix:verification}), and, later, the conclusions (\S
\ref{sec:results-on-example}).

\begin{figure}
\resizebox{\textwidth}{!}{\includegraphics{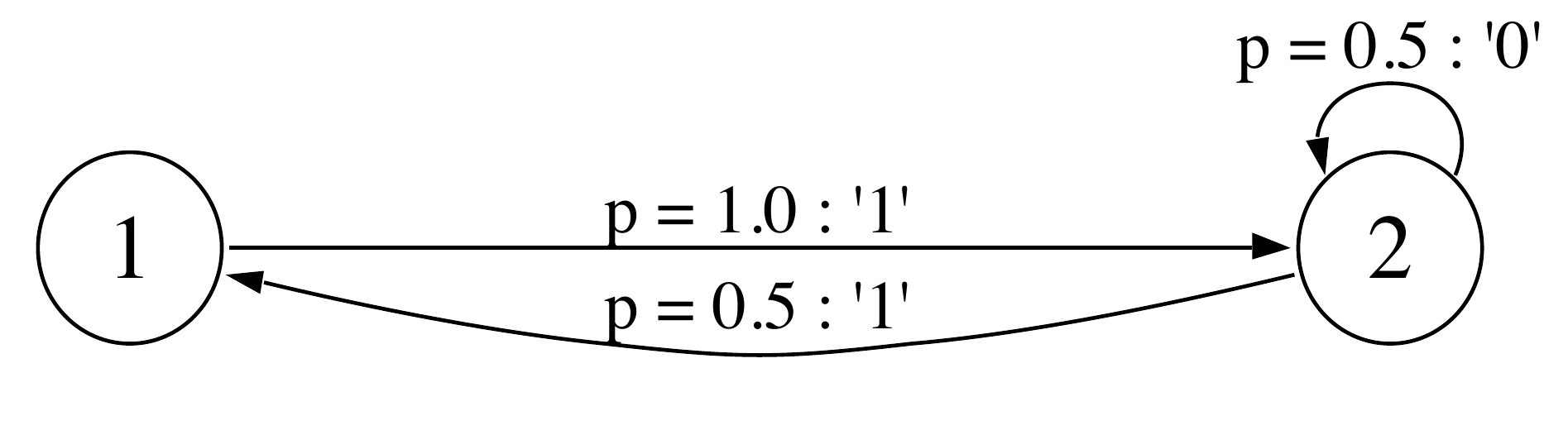}}
\caption{State-transition diagram for the ``even process''.  The legends on the
  transition arrows indicate the probability of making the transition, and the
  observation which occurs when the transition happens.  The observation
  $X_t=1$ when entering or leaving state 1, otherwise it is 0.  This creates
  blocks of 1s of even length, separated by blocks of 0s of arbitrary length.
  The result is a finite-state process which is not a Markov chain of any
  order.}
\label{fig:even-process}
\end{figure}

The data-generating process $\TrueMeasure$ is a stationary and ergodic measure
on the space of binary sequences, i.e., $\ObservableSpace = \left\{ 0,
  1\right\}$, and the $\sigma$-field $\ObservableSigmaField =
2^\ObservableSpace$.  The measure is naturally represented as a function of a
two-state Markov chain $S_1^{\infty}$, with $S_t \in \left\{1,2\right\}$.  The
transition matrix is
\[
\mathbf{T} = \left[\begin{array}{cc} 0.0 & 1.0 \\
0.5 & 0.5 \end{array}\right]
\]
so that the invariant distribution puts probability $1/3$ on state 1 and
probability $2/3$ on state 2; take $S_1$ to be distributed accordingly.  The
observed process is a binary function of the latent state transitions, $X_t =
0$ if $S_t = S_{t+1} = 2$ and $X_t = 1$ otherwise.  Figure
\ref{fig:even-process} depicts the transition and observation structure.
Qualitatively, $X_1^{\infty}$ consists of blocks of $1$s of even length,
separated by blocks of $0$s of arbitrary length.  Since the joint process
$\left\{(S_t,X_t)\right\}_{1\leq t\leq\infty}$ is a stationary and ergodic
Markov chain, $X_1^{\infty}$ is also stationary, ergodic and mixing.

This stochastic process comes from symbolic dynamics
\citep{Lind-Marcus,Kitchens}, where it is known as the ``even process'', and
serves as a basic example of the class of {\em sofic} processes
\citep{Weiss-1973}, which have finite Markovian representations, as in Figure
\ref{fig:even-process}, but are not Markov at any finite order.  (If $X_t =1,
X_{t-1} = 1, \ldots X_{t-k}=1$ for any finite $k$, the corresponding $S_{t-i}$
must have alternated between one and two, but whether $S_t$ is one or two, and
thus the distribution of $X_{t+1}$, cannot be determined from the length-$k$
history alone.)  More exactly \citep{Kitchens-Tuncel}, sofic systems or
``finitary measures'' are ones which are images of Markov chains under factor
maps, and strictly sofic systems, such as the even process, are sofic systems
which are not themselves Markov chains of any order.  Despite their simplicity,
these models arise naturally when studying the time series of chaotic dynamical
systems
\citep{Badii-Politi,JPC-semantics,CMPPSS,Daw-Finney-Tracy-symbolic-dynamics},
as well as problems in statistical mechanics
\citep{Perry-Binder-finite-stat-compl} and crystallography
\citep{Varn-infinite-range-order}.

Let $\ParameterSpace_k$ be the space of all binary Markov chains of order $k$
with strictly positive transition probabilities and their respective stationary
distributions; each $\ParameterSpace_k$ has dimension $2^k$.  (Allowing some
transition probabilities to be zero creates uninteresting technical
difficulties.)  Since each hypothesis is equivalent to a function
$\ObservableSpace^{k+1} \mapsto (0,1]$, we can give $\ParameterSpace_k$ the
topology of pointwise convergence of functions, and the corresponding Borel
$\sigma$-field.  We will take $\ParameterSpace = \bigcup_{k =
  1}^{\infty}{\ParameterSpace_k}$, identifying $\ParameterSpace_k$ with the
appropriate subset of $\ParameterSpace_{k+1}$.  Thus $\ParameterSpace$ consists
of all strictly-positive stationary binary Markov chains, of whatever order,
and is infinite-dimensional.

As for the prior $\Posterior_0$, it will be specified in more detail below (\S
\ref{sec:verification-of-assumptions}).  At the very least, however, it needs
to have the ``Kullback-Leibler rate property'', i.e., to give positive
probability to every $\epsilon$ ``neighborhood''
$N_{\epsilon}(\parametervalue)$ around every $\parametervalue \in
\ParameterSpace$, i.e., the set of hypotheses whose Kullback-Leibler divergence
from $\parametervalue$ grows no faster than $\epsilon$:
\[
N_{\epsilon}(\parametervalue) =\left\{ \parametervalue^{\prime}: \epsilon \geq \lim_{t\rightarrow\infty}{\frac{1}{t}\int{dx_1^t \modeldensity(x_1^t) \log{\frac{\modeldensity(x_1^t)}{\altmodeldensity(x_1^t)}} }} \right\}
\]
(The limit exists for all $\parametervalue,\parametervalue^{\prime}$
combinations \citep{Gray-entropy}.)

This example is simple, but it is also beyond the scope of existing work on
Bayesian convergence in several ways.  First, the data-generating process
$\TrueMeasure$ is not even Markov.  Second, $\TrueMeasure \not\in
\ParameterSpace$, so all the hypotheses are wrong, and the truth is certainly
not in the support of the prior.  ($\TrueMeasure$ can however be approximated
arbitrarily closely, in various process metrics, by distributions from
$\ParameterSpace$.)  Third, because $\TrueMeasure$ is ergodic, and ergodic
distributions are extreme points in the space of stationary distributions
\citep{Dynkin-suff-stats-and-extreme-points}, it cannot be represented as a
mixture of distributions in $\ParameterSpace$.  This means that the Doob-style
theorem of Ref.\ \citep{Bayesian-consistency-for-stationary-models} does not
apply, and even the subjective certainty of convergence is not assured.  The
results of Refs.\
\citep{Kleijn-van-der-Vaart,Zhang-from-epsilon-to-KL,Berk-limiting-behavior-of-posterior,Berk-consistency}
on mis-specified models do not hold because the data are dependent. To be as
concrete and explicit as possible, the analysis here will focus on the even
process, but only the constants would change if $\TrueMeasure$ were any other
strictly sofic process.  Much of it would apply even if $\TrueMeasure$ were a
stochastic context-free language or pushdown automaton
\citep{Charniak-statistical-language-learning}, where in effect the number of
hidden states is infinite, though some of the details in Appendix
\ref{appendix:verification} would change.

Ref.\ \citep{Ornstein-Weiss-how-sampling-reveals-a-process} describes a
non-parametric procedure which will adaptively learn to predict a class of
discrete stochastic processes which includes the even process.  Ref.\
\citep{CSSR-for-UAI} introduces a frequentist algorithm which consistently
reconstructs the hidden-state representation of sofic processes, including the
even process.  Ref.\ \citep{Strelioff-JPC-Hubler-Markov-Bayes} considers
Bayesian estimation of the even process, using Dirichlet priors for
finite-order Markov chains, and employing Bayes factors to decide which order
of chain to use for prediction.

\subsection{Assumptions}

The needed assumptions have to do with the dynamical properties of the data
generating process $\TrueMeasure$, and with how well the dynamics meshes both
with the class of hypotheses $\ParameterSpace$ and with the prior distribution
$\Posterior_0$ over those hypotheses.

\begin{assumption}
  The likelihood ratio $\LikelihoodRatio_t(\parametervalue)$ is $\Filtration{t}
  \times \ParameterSigmaField$-measurable for all $t$.
  \label{assumption:measurable-likelihood}
\end{assumption}
%assumption:1

The next two assumptions actually need only hold for $\Posterior_0$-almost-all
$\parametervalue$.  But this adds more measure-0 caveats to the results, and it
is hard to find a natural example where it would help.

\begin{assumption}
  For every $\parametervalue \in \ParameterSpace$, the Kullback-Leibler
  divergence rate from $\TrueMeasure$,
  \[
  \divrate(\parametervalue) =
  \lim_{t\rightarrow\infty}{\frac{1}{t}\Expect{\log{\frac{\truedensity(X_1^t)}{\modeldensity(X_1^t)}}}}
  \]
  exists (possibly being infinite) and is
  $\ParameterSigmaField$-measurable.
  \label{assumption:divergence-rates-exist}
\end{assumption}
%assumption:2

As mentioned, any two distinct ergodic measures are mutually singular, so there
is a consistent test which can separate them.
(\citep{Ryabko-Ryabko-testing-ergodic} constructs an explicit but not
necessarily optimal test.)  One interpretation of the divergence rate
\citep{Gray-entropy} is that it measures the maximum exponential rate at which
the power of such tests approaches 1, with $d=0$ and $d=\infty$ indicating sub-
and supra- exponential convergence, respectively.

\begin{assumption}
  For each $\parametervalue \in \ParameterSpace$, the generalized or relative
  asymptotic equipartition property holds, and so
  \begin{equation}
    \lim{\frac{1}{t}\log{\LikelihoodRatio_t(\parametervalue)}} = -\divrate(\parametervalue)
    \label{eqn:relative-AEP}
  \end{equation}
  with $\TrueMeasure$-probability 1.
  \label{assumption:relative-AEP-holds}
\end{assumption}
%assumption:3

Refs.\ \citep{Algoet-and-Cover-on-AEP,Gray-entropy} give sufficient, but not
necessary, conditions sufficient for Assumption
\ref{assumption:relative-AEP-holds} to hold for a given $\parametervalue$.  The
ordinary, non-relative asymptotic equipartition property, also known as the
Shannon-McMillan-Breiman theorem, is that
$\lim{t^{-1}\log{\truedensity(x_1^t)}} = -h_{\TrueMeasure}$ a.s., where
$h_{\TrueMeasure}$ is the entropy rate of the data-generating process.  (See
\citep{Gray-entropy}.)  If this holds and $h_{\TrueMeasure}$ is finite, one
could rephrase Assumption \ref{assumption:relative-AEP-holds} as
$\lim{t^{-1}\log{\modeldensity(X_1^t)}} =
-h_{\TrueMeasure}-\divrate(\parametervalue)$ a.s., and state results in terms
of the likelihood rather than the likelihood ratio.  (Cf.\ \cite[ch.\
5]{Andy-Fraser-on-HMMs}.)  However, there are otherwise-well-behaved processes
for which $h_{\TrueMeasure}=-\infty$, at least in the usual choice of reference
measure, so I will restrict myself to likelihood ratios.

The meaning of Assumption \ref{assumption:relative-AEP-holds} is that, relative
to the true distribution, the likelihood of each $\parametervalue$ goes to zero 
exponentially, the rate being the Kullback-Leibler divergence rate.  Roughly
speaking, an integral of exponentially-shrinking quantities will tend to be
dominated by the integrand with the slowest rate of decay.  This suggests that
the posterior probability of a set $A \subseteq \ParameterSpace$ depends on the
smallest divergence rate which can be attained at a point of prior support
within $A$.  Thus, adapting notation from large deviations theory, define
\begin{eqnarray*}
\divrate(A) & \equiv & \essinf_{\parametervalue \in A}{\divrate(\parametervalue)}\\
\Action(\parametervalue) & \equiv & \divrate(\parametervalue) - \divrate(\ParameterSpace)\\
\Action(A) & \equiv & \essinf_{\parametervalue \in A}{\Action(\parametervalue)}
\end{eqnarray*}
where here and throughout $\essinf$ is the essential infimum with respect to
$\Posterior_0$, i.e., the greatest lower bound which holds with
$\Posterior_0$-probability 1.

Our further assumptions are those needed for the ``roughly speaking'' and
``should'' statements of the previous paragraph to be true, so that, for
reasonable sets $A \in \ParameterSigmaField$,
\[
\lim{\frac{1}{t}\log{\prioravg{\LikelihoodRatio_t A}}} = -\divrate(A)
\]

Let $\InfinitySet \equiv \left\{\parametervalue:~ \divrate(\parametervalue) = \infty\right\}$.
\begin{assumption}
  $\Posterior_0(\InfinitySet) < 1$
  \label{assumption:best-divergence-is-finite}
\end{assumption}
%assumption:4

If this assumption fails, then every hypothesis in the support of the prior
doesn't just diverge from the true data-generating distribution, it diverges so
rapidly that the error rate of a test against the latter distribution goes to
zero faster than any exponential.  (One way this can happen is if every
hypothesis has a finite-dimensional distribution assigning probability zero to
some event of positive $\TrueMeasure$-probability.)  The methods of this paper
seem to be of no use in the face of such extreme mis-specification.

Our first substantial assumption is that the prior distribution does not give
too much weight to parts of $\ParameterSpace$ where the log likelihood
converges badly.

\begin{assumption}
  There exists a sequence of sets $G_t \rightarrow \ParameterSpace$ such that
  \begin{enumerate}
  \item $\Posterior_0(G_t) \geq 1 - \alpha \myexp{- t \beta}$, for some $\alpha > 0$,
    $\beta > 2 \divrate(\ParameterSpace)$;
  \item The convergence of Eq.\ \ref{eqn:relative-AEP} is uniform in
    $\parametervalue$ over $G_t \setminus \InfinitySet$;
  \item $\divrate(G_t) \rightarrow \divrate(\ParameterSpace)$.
  \end{enumerate}
  \label{assumption:exponentially-egorov}
\end{assumption}
%assumption:5

{\em Comment 1:} An analogy with the method of sieves
\citep{Geman-and-Hwang-on-method-of-sieves} may clarify the meaning of the
assumption.  If we were constrained to some fixed $G$, the uniform convergence
in the second part of the assumption would make the convergence of the
posterior distribution fairly straightforward.  Now imagine that the constraint
set is gradually relaxed, so that at time $t$ the posterior is confined to
$G_t$, which grows so slowly that convergence is preserved.  (Assumption
\ref{assumption:sufficiently-rapid-convergence} below is, in essence, about the
relaxation being sufficiently slow.)  The theorems work by showing that the
behavior of the posterior distribution on the full space $\ParameterSpace$ is
dominated by its behavior on this ``sieve''.

{\em Comment 2:} Recall that by Egorov's theorem \cite[Lemma 1.36, p.\
18]{Kallenberg-mod-prob}, if a sequence of finite, measurable functions
$f_t(\parametervalue)$ converges pointwise to a finite, measurable function
$f(\parametervalue)$ for $\Posterior_0$-almost-all $\parametervalue \in G$,
then for each $\epsilon > 0$, there is a (possibly empty) $B \subset G$ such
that $\Posterior_0(B) \leq \epsilon$, and the convergence is uniform on $G
\setminus B$.  Thus the first two parts of the assumption really follow for
free from the measurability in $\parametervalue$ of likelihoods and divergence
rates.  (That $\beta$ needs to be at least $2\divrate(\ParameterSpace)$ becomes
apparent in the proof of Lemma
\ref{lemma:integrated-likelihood-mostly-from-good-set}, but that could always
be arranged.)  The extra content comes in the third part of the assumption,
which could fail if the lowest-divergence hypotheses were also the ones where
the convergence was slowest, consistently falling into the bad sets $B$ allowed
by Egorov's theorem.

For each measurable $A \subseteq \ParameterSpace$, for every $\delta > 0$,
there exists a random natural number $\convtime(A,\delta)$ such that
\[
t^{-1}\log{\prioravg{A \LikelihoodRatio_t}} \leq \delta + \limsup_{t}{t^{-1}\log{\prioravg{A \LikelihoodRatio_t}}}
\]
for all $t > \convtime(A,\delta)$, provided the $\limsup$ is finite.  We need
this random last-entry time $\convtime(A,\delta)$ to state the next assumption.

\begin{assumption}
  The sets $G_t$ of the previous assumption can be chosen so that, for every
  $\delta$, the inequality $t \geq \convtime(G_t, \delta)$ holds a.s. for all
  sufficiently large $t$.
  \label{assumption:sufficiently-rapid-convergence}
\end{assumption}
%assumption:6
The fraction of the prior probability mass outside of $G_t$ is exponentially
small in $t$, with the decay rate large enough that (Lemma
\ref{lemma:integrated-likelihood-mostly-from-good-set}) the posterior
probability mass outside $G_t$ also goes to zero.  Using the analogy to the
sieve again, the meaning of the assumption is that the convergence of the
log-likelihood ratio is sufficiently fast, and the relaxation of the sieve is
sufficiently slow, that, at least eventually, every set $G_t$ has
$\delta$-converged by $t$, the time when we start using it.

To show convergence of the posterior measure, we need to be able to control the
convergence of the log-likelihood on sets smaller than the whole parameter
space.
\begin{assumption}
  The sets $G_t$ of the previous two assumptions can be chosen so that, for any
  set $A$ with $\Posterior_0(A) > 0$, $\divrate(G_t \cap A) \rightarrow
  \divrate(A)$.
  \label{assumption:good-sets-are-good-even-for-subsets}
\end{assumption}
%assumption:7

Assumption \ref{assumption:good-sets-are-good-even-for-subsets} could be
replaced by the logically-weaker assumption that for each set $A$, there exist
a sequence of sets $G_{t,A}$ satisfying the equivalents of Assumptions
\ref{assumption:exponentially-egorov} and
\ref{assumption:sufficiently-rapid-convergence} for the prior measure
restricted to $A$.  Since the most straight-forward way to check such an
assumption would be to verify Assumption
\ref{assumption:good-sets-are-good-even-for-subsets} as stated, the extra
generality does not seem worth it.

\subsubsection{Verification of Assumptions for the Example}
\label{sec:verification-of-assumptions}

Since every $\parametervalue \in \ParameterSpace$ is a finite-order Markov
chain, and $\TrueMeasure$ is stationary and ergodic, Assumption
\ref{assumption:measurable-likelihood} is unproblematic, while Assumptions
\ref{assumption:divergence-rates-exist} and \ref{assumption:relative-AEP-holds}
hold by virtue of \citep{Algoet-and-Cover-on-AEP}.

It is easy to check that $\inf_{\parametervalue \in
  \ParameterSpace_k}{\divrate(\parametervalue)} > 0$ for each $k$.  (The
infimum is not in general attained by any $\parametervalue \in
\ParameterSpace_k$, though it could be if the chains were allowed to have some
transition probabilities equal to zero.)  The infimum over $\ParameterSpace$ as
a whole, however, is zero.  Also, $\divrate(\parametervalue) < \infty$
everywhere (because none of the hypotheses' transition probabilities are zero),
so the possible set $I$ of $\parametervalue$ with infinite divergence rates is
empty, disposing of Assumption \ref{assumption:best-divergence-is-finite}.

Verifying the remaining assumptions means building a sequence $G_t$ of
increasing subsets of $\ParameterSpace$ on which the convergence of
$t^{-1}\log{\LikelihoodRatio_t}$ is uniform and sufficiently rapid, and
ensuring that the prior probability of these sets grows fast enough.  This will
be done by exploiting some finite-sample deviation bounds for the even process,
which in turn rest on its mixing properties and the method of types.  Details
are referred to Appendix \ref{appendix:verification}.  The upshot is that the
sets $G_t$ consist of chains whose order is less than or equal to
$\frac{\log{t}}{2/3 + \epsilon} - 1$, for some $\epsilon > 0$, and where the
absolute logarithm of all the transition probabilities is bounded by $C
t^{\gamma}$, where the positive constant $C$ is arbitrary but $0 < \gamma <
\frac{2/3+\epsilon/2}{2/3 +\epsilon}$.  (With a different strictly sofic
process $\TrueMeasure$, the constant $2/3$ in the preceding expressions should
be replaced by $h_{\TrueMeasure}$.)  The exponential rate $\beta > 0$ for the
prior probability of $G_t^c$ can be chosen to be arbitrarily small.

\section{Results}
\label{sec:asymptotics}

I first give the theorems here, without proof.  The proofs, in \S\S
\ref{sec:without-prior-control}--\ref{sec:rate-of-convergence}, are accompanied
by re-statements of the theorems, for the reader's convenience.

There are six theorems.  The first upper-bounds the growth rate of the
posterior density at a given point $\parametervalue$ in $\ParameterSpace$.  The
second matches the upper bound on the posterior density with a lower bound,
together providing the growth-rate for the posterior density.  The third is
that $\Posterior_t(A) \rightarrow 0$ for any set with $J(A) > 0$, showing that
the posterior concentrates on the divergence-minimizing part of the hypothesis
space.  The fourth is a kind of large deviations principle for the posterior
measure.  The fifth bounds the asymptotic Hellinger and total variation
distances between the posterior predictive distribution and the actual
conditional distribution of the next observation.  Finally, the sixth theorem
establishes rates of convergence.

The first result uses only Assumptions
\ref{assumption:measurable-likelihood}--\ref{assumption:relative-AEP-holds}.
(It is not very interesting, however, unless
\ref{assumption:best-divergence-is-finite} is also true.)  The latter three,
however, all depend on finer control of the integrated likelihood, and so finer
control of the prior, as embodied in Assumptions
\ref{assumption:exponentially-egorov}--\ref{assumption:sufficiently-rapid-convergence}.
More exactly, those additional assumptions concern the interplay between the
prior and the data-generating process, restricting the amount of prior
probability which can be given to hypotheses whose log-likelihoods converge
excessively slowly under $\TrueMeasure$.  I build to the first result in the
next sub-section, then turn to the control of the integrated likelihood and its
consequences in the next three sub-sections, and then consider how these
results apply to the example.

\begin{theorem}
  Under Assumptions
  \ref{assumption:measurable-likelihood}--\ref{assumption:relative-AEP-holds},
  with probability 1, for all $\parametervalue$ where
  $\posteriordensity_0(\parametervalue) > 0$,
  \[
    \limsup_{t\rightarrow\infty}{\frac{1}{t}\log{\posteriordensity_t(\parametervalue)}} \leq -\Action(\parametervalue)
  \]
\end{theorem}

\begin{theorem}
  Under Assumptions
  \ref{assumption:measurable-likelihood}--\ref{assumption:sufficiently-rapid-convergence},
  for all $\parametervalue \in \ParameterSpace$ where
  $\posteriordensity_0(\parametervalue) > 0$,
  \[
  \lim_{t\rightarrow\infty}{\frac{1}{t}\log{\posteriordensity_t(\parametervalue)}}
  = -\Action(\parametervalue)
  \]
  with probability 1.
\end{theorem}

\begin{theorem}
  Make Assumptions
  \ref{assumption:measurable-likelihood}--\ref{assumption:good-sets-are-good-even-for-subsets}.
  Pick any set $A \in \ParameterSigmaField$ where $\Posterior_0(A) > 0$ and
  $\divrate(A) > \divrate(\ParameterSpace)$.  Then $\Posterior_t(A) \rightarrow
  0$ a.s.
\end{theorem}

\begin{theorem}
  Under the conditions of Theorem \ref{thm:convergence-of-measure}, if
  $A\in\ParameterSigmaField$ is such that
  \[
  -\limsup{t^{-1}\log{\Posterior_0(A \cap G_t^c)}} = \beta^{\prime} \geq
  2\divrate(A)
  \]
  then
  \[
  \lim_{t\rightarrow\infty}{\frac{1}{t}\log{\Posterior_t(A)}} =
  -\Action(A)
  \]
  In particular, this holds whenever $2\divrate(A) < \beta$ or $A \subset
  \bigcap_{k=n}^{\infty}{G_k}$ for some $n$.
\end{theorem}

\begin{theorem}
  Under Assumptions
  \ref{assumption:measurable-likelihood}--\ref{assumption:good-sets-are-good-even-for-subsets}, with probability 1,
  \begin{eqnarray*}
    \limsup_{t\rightarrow\infty}{\rho^2_{H}(\TrueMeasure^t,\PredictiveMeasure^t)} & \leq & \divrate(\ParameterSpace)\\
    \limsup_{t\rightarrow\infty}{\rho_{TV}^2(\TrueMeasure^t,\PredictiveMeasure^t)} & \leq & 4\divrate(\ParameterSpace)
  \end{eqnarray*}
  where $\rho_{H}$ and $\rho_{TV}$ are, respectively, the Hellinger and total
  variation metrics.
\end{theorem}

\begin{theorem}
  Make assumptions
  \ref{assumption:measurable-likelihood}--\ref{assumption:good-sets-are-good-even-for-subsets},
  and pick a positive sequence $\epsilon_t$ where $\epsilon_t \rightarrow 0$,
  $t\epsilon_t \rightarrow\infty$.  If, for each $\delta > 0$, 
  \[
  \tau(G_t \cap N^c_{\epsilon_t},\delta) \leq t
  \]
  eventually almost surely, then
  \[
    \Posterior_t(N_{\epsilon_t}) \rightarrow 1
  \]
  with probability 1.
\end{theorem}
%theorem:6

\setcounter{theorem}{0}

\subsection{Upper Bound on the Posterior Density}
\label{sec:without-prior-control}

The primary result of this section is a pointwise upper bound on the growth
rate of the posterior density.  To establish it, I use some subsidiary lemmas,
which also recur in later proofs.  Lemma
\ref{lemma:convergence-for-almost-all-parameters} extends the almost-sure
convergence of the likelihood (Assumption \ref{assumption:relative-AEP-holds})
from holding pointwise in $\ParameterSpace$ to holding simultaneously for all
$\parametervalue$ on a (possibly random) set of $\Posterior_0$-measure 1.  Lemma
\ref{lemma:integrated-likelihood-tracks-the-best} shows that the prior-weighted
likelihood ratio, $\prioravg{\LikelihoodRatio_t}$ tends to be at least
$\myexp{-t\divrate(\ParameterSpace)}$.  (Both assertions are made more precise
in the lemmas themselves.)

I begin with a proposition about exchanging the order of universal quantifiers
(with almost-sure caveats).

\begin{lemma}
  Let $Q \subset \ParameterSpace \times \Omega$ be jointly measurable, with
  sections $Q_\parametervalue = \left\{ \omega:~ (\omega,\parametervalue) \in
    Q\right\}$ and $Q_\omega = \left\{\parametervalue:~
    (\omega,\parametervalue) \in Q\right\}$.  If, for some probability measure
  $\TrueMeasure$ on $\Omega$,
  \begin{equation}
    \forall \parametervalue \TrueMeasure\left(Q_{\parametervalue}\right) = 1
  \end{equation}
  then for any probability measure
  $\Posterior$ on $\ParameterSpace$
  \begin{equation}
    \TrueMeasure\left(\left\{\omega: \Posterior\left(Q_{\omega}\right) = 1\right\}\right) = 1
  \end{equation}
  \label{lemma:quantifier-reversal}
\end{lemma}
%lemma:1

In words, if, for all $\parametervalue$, some property holds a.s., then
a.s. the property holds {\em simultaneously} for almost all $\parametervalue$.

\textsc{Proof:} Since $Q$ is measurable, for all $\omega$ and
$\parametervalue$, the sections are measurable, and the measures of the
sections, $\TrueMeasure(Q_{\parametervalue})$ and $\Posterior(Q_\omega)$, are
measurable functions of $\parametervalue$ and $\omega$, respectively.  Using
Fubini's theorem,
\begin{eqnarray*}
\int_{\ParameterSpace}{\TrueMeasure(Q_{\parametervalue}) d\Posterior(\parametervalue)} 
& = & \int_{\ParameterSpace}{\int_{\Omega}{\mathbf{1}_{Q}(\omega,\parametervalue) d\TrueMeasure(\omega)}d\Posterior(\parametervalue)}\\
& = & \int_{\Omega}{\int_{\ParameterSpace}{\mathbf{1}_{Q}(\omega,\parametervalue) d\Posterior(\parametervalue)}d\TrueMeasure(\omega)}\\
& = & \int_{\Omega}{\Posterior(Q_\omega) d\TrueMeasure(\omega)}
\end{eqnarray*}
By hypothesis, however, $\TrueMeasure(Q_{\parametervalue}) = 1$ for all
$\parametervalue$.  Hence it must be the case that $\Posterior(Q_\omega) =1$
for $\TrueMeasure$-almost-all $\omega$.  (In fact, the set of $\omega$ for
which this is true must be a measurable set.) $\Box$

\begin{lemma}
  Under Assumptions
  \ref{assumption:measurable-likelihood}--\ref{assumption:relative-AEP-holds},
  there exists a set $C \subseteq \ObservableSpace^{\infty}$, with
  $\TrueMeasure(C) = 1$, where, for every $y \in C$, there exists a $Q_y \in
  \ParameterSigmaField$ such that, for every $\parametervalue \in Q_y$, Eq.\
  \ref{eqn:relative-AEP} holds.  Moreover, $\Posterior_0(Q_y) =1$.
  \label{lemma:convergence-for-almost-all-parameters}
\end{lemma}
%lemma:2

\textsc{Proof}: Let the set $Q$ consist of the $\parametervalue,\omega$ pairs
where Eq.\ \ref{eqn:relative-AEP} holds, i.e., for which
\[
\lim{\frac{1}{t}\log{\LikelihoodRatio_t(\parametervalue,\omega)}} = -\divrate(\parametervalue)~,
\]
being explicit about the dependence of the likelihood ratio on $\omega$.
Assumption \ref{assumption:relative-AEP-holds} states that
$\forall\parametervalue \TrueMeasure(Q_\parametervalue) = 1$, so applying Lemma
\ref{lemma:quantifier-reversal} just needs the verification that $Q$ is jointly
measurable.  But, by Assumptions \ref{assumption:measurable-likelihood} and
\ref{assumption:divergence-rates-exist}, $\divrate(\cdot)$ is
$\ParameterSigmaField$-measurable, and $\LikelihoodRatio_t(\parametervalue)$ is
$\Filtration{t}\times\ParameterSigmaField$-measurable for each $t$, so the set
$Q$ where the convergence holds are
$\Filtration{\infty}\times\ParameterSigmaField$-measurable.  Everything then
follows from the preceding lemma.  $\Box$

{\em Remark:} Lemma \ref{lemma:convergence-for-almost-all-parameters}
generalizes Lemma 3 in \citep{Barron-Schervish-Wasserman}.  Lemma
\ref{lemma:quantifier-reversal} is a specialization of the quantifier-reversal
lemma used in \citep{McAllister-PAC-Bayes} to prove PAC-Bayesian theorems for
learning classifiers.  Lemma \ref{lemma:quantifier-reversal} could be used to
extend any of the results below which hold a.s. for each $\parametervalue$ to
ones which a.s. hold {\em simultaneously} almost everywhere in
$\ParameterSpace$.  This may seem too good to be true, like an alchemist's
recipe for turning the lead of pointwise limits into the gold of uniform
convergence.  Fortunately or not, however, the lemma tells us nothing about the
{\em rate} of convergence, and is compatible with its varying across
$\ParameterSpace$ from instantaneous to arbitrarily slow, so uniform laws need
stronger assumptions.

\begin{lemma}
  Under Assumptions
  \ref{assumption:measurable-likelihood}--\ref{assumption:relative-AEP-holds},
  for every $\epsilon > 0$, it is almost sure that the ratio between the
  integrated likelihood and the true probability density falls below
  $\myexp{-t(\divrate(\ParameterSpace) + \epsilon)}$ only finitely often:
  \begin{equation}
    \TrueMeasure\left\{x_1^\infty:~ \prioravg{\LikelihoodRatio_t} \leq \myexp{-t(\divrate(\ParameterSpace) + \epsilon)}, ~\mathrm{i.o.}\right\} = 0
    \label{eqn:lower-bound-on-integrated-likelihood}
  \end{equation}
  and as a corollary, with probability 1,
  \begin{equation}
    \liminf_{t\rightarrow\infty}{\frac{1}{t}\log{\prioravg{\LikelihoodRatio_t}}} \geq -\divrate(\ParameterSpace)
  \end{equation}
  \label{lemma:integrated-likelihood-tracks-the-best}
\end{lemma}
%lemma:3

\textsc{Proof:} It's enough to show that Eq.\
\ref{eqn:lower-bound-on-integrated-likelihood} holds for all $x_1^\infty$ in
the set $C$ from the previous lemma, since that set has probability 1.

Let $N_{\epsilon/2}$ be the set of all $\parametervalue$ in the support of
$\Posterior_0$ such that $\divrate(\parametervalue) \leq
\divrate(\ParameterSpace) + \epsilon/2$.  Since $x_1^\infty \in C$, the
previous lemma tells us there exists a set $Q_{x_1^\infty}$ of
$\parametervalue$ for which Eq.\ \ref{eqn:relative-AEP} holds under the
sequence $x_1^\infty$.

\begin{eqnarray*}
\lefteqn{\myexp{t(\epsilon+\divrate(\ParameterSpace))} \prioravg{\LikelihoodRatio_t} =  \int_{\ParameterSpace}{\LikelihoodRatio_t(\parametervalue)\myexp{t(\epsilon+\divrate(\ParameterSpace))} d\Posterior_0(\parametervalue)}} & &\\
& \geq & \int_{N_{\epsilon/2} \cap Q_{x_1^{\infty}}}{\LikelihoodRatio_t(\parametervalue)\myexp{t(\epsilon+\divrate(\ParameterSpace))} d\Posterior_0(\parametervalue)}\\
& = & \int_{N_{\epsilon/2} \cap Q_{x_1^{\infty}}}{\myexp{t\left[\epsilon + \divrate(\ParameterSpace) + \frac{\log{\LikelihoodRatio_t(\parametervalue)}}{t}\right]}d\Posterior_0(\parametervalue)}
\end{eqnarray*}
By Assumption \ref{assumption:relative-AEP-holds},
\[
\lim_{t\rightarrow\infty}{\frac{1}{t}\log{\LikelihoodRatio_t(\parametervalue)} } = -\divrate(\parametervalue)
\]
and for all $\parametervalue \in N_{\epsilon/2}$, $\divrate(\parametervalue)
\leq \divrate(\ParameterSpace) + \epsilon/2$, so
\[
\liminf_{t\rightarrow\infty}{\myexp{t\left[\epsilon + \divrate(\ParameterSpace) + \frac{1}{t}\log{\LikelihoodRatio_t(\parametervalue)} \right]}} = \infty
\]
a.s., for all $\parametervalue \in N_{\epsilon/2} \cap Q_{x_1^{\infty}}$.  We
must have $\Posterior_0(N_{\epsilon/2}) > 0$, otherwise
$\divrate(\ParameterSpace)$ would not be the essential infimum, and we know
from the previous lemma that $\Posterior_0(Q_{x_1^\infty}) = 1$, so
$\Posterior_0(N_{\epsilon/2} \cap Q_{x_1^{\infty}}) > 0$.  Thus, Fatou's lemma
gives
\[
\lim_{t\rightarrow\infty}{\int_{N_{\epsilon/2} \cap Q_{x_1^{\infty}}}{\myexp{t\left[\epsilon + \divrate(\ParameterSpace) + \frac{1}{t}\log{\LikelihoodRatio_t(\parametervalue)}\right]}d\Posterior_0(\parametervalue)}} = \infty
\]
so
\[
\lim_{t\rightarrow\infty}{\myexp{t(\epsilon+\divrate(\ParameterSpace))}\prioravg{\LikelihoodRatio_t}} = \infty
\]
and hence
\begin{equation}
\prioravg{\LikelihoodRatio_t} > \myexp{-t(\epsilon+\divrate(\ParameterSpace))}
\label{eqn:integrated-likelihood-tracks-true-density}
\end{equation}
for all but finitely many $t$.  Since this holds for all $x_1^\infty \in C$,
and $\TrueMeasure(C) = 1$, Equation
\ref{eqn:integrated-likelihood-tracks-true-density} holds a.s., as was to be
shown. The corollary statement follows immediately. $\Box$

\begin{theorem}
  Under Assumptions
  \ref{assumption:measurable-likelihood}--\ref{assumption:relative-AEP-holds},
  with probability 1, for all $\parametervalue$ where
  $\posteriordensity_0(\parametervalue) > 0$,
  \begin{equation}
    \limsup_{t\rightarrow\infty}{\frac{1}{t}\log{\posteriordensity_t(\parametervalue)}} \leq -\Action(\parametervalue)
  \end{equation}
  \label{thm:upper-bound-on-long-run-fitness}
\end{theorem}
%theorem:1

\textsc{Proof:} As remarked,
\[
\posteriordensity_t(\parametervalue) = \posteriordensity_0(\parametervalue) \frac{\LikelihoodRatio_t(\parametervalue)}{\prioravg{\LikelihoodRatio_t}}
\]
so
\[
\frac{1}{t}\log{\posteriordensity_t(\parametervalue)} = \frac{1}{t}\log{\posteriordensity_0(\parametervalue)} + \frac{1}{t}\log{\LikelihoodRatio_t(\parametervalue)} - \frac{1}{t}\log{\prioravg{\LikelihoodRatio_t}}
\]

By Assumption \ref{assumption:relative-AEP-holds}, for each $\epsilon > 0$,
it's almost sure that
\[
\frac{1}{t}\log{\LikelihoodRatio_t(\parametervalue)} \leq -\divrate(\parametervalue) + \epsilon/2
\]
for all sufficiently large $t$, while by Lemma
\ref{lemma:integrated-likelihood-tracks-the-best}, it's almost sure that
\[
\frac{1}{t}\log{\prioravg{\LikelihoodRatio_t}} \geq -\divrate(\ParameterSpace) - \epsilon/2
\]
for all sufficiently large $t$.  Hence, with probability 1,
\[
\frac{1}{t}\log{\posteriordensity_t(\parametervalue)} \leq \divrate(\ParameterSpace)  - \divrate(\parametervalue) +\epsilon + \frac{1}{t}\log{\posteriordensity_0(\parametervalue)}
\]
for all sufficiently large $t$.  Hence
\[
\limsup_{t\rightarrow\infty}{\frac{1}{t}\log{\posteriordensity_t(\parametervalue)}} \leq \divrate(\ParameterSpace) - \divrate(\parametervalue) = -\Action(\parametervalue)
\]
$\Box$

Lemma \ref{lemma:integrated-likelihood-tracks-the-best} gives a lower bound on
the integrated likelihood ratio, showing that in the long run it has to be at
least as big as $\myexp{-t\divrate(\ParameterSpace))}$.  (More precisely, it is
significantly smaller than that on vanishingly few occasions.)  It does not,
however, rule out being larger.  Ideally, we would be able to match this lower
bound with an upper bound of the same form, since $\divrate(\ParameterSpace)$
is the best attainable divergence rate, and, by Lemma
\ref{lemma:convergence-for-almost-all-parameters}, log likelihood ratios per
unit time are converging to divergence rates for $\Posterior_0$-almost-all
$\parametervalue$, so values of $\parametervalue$ for which
$\divrate(\parametervalue)$ are close to $\divrate(\ParameterSpace)$ should
come to dominate the integral in $\prioravg{\LikelihoodRatio_t}$.  It would
then be fairly straightforward to show convergence of the posterior
distribution.

Unfortunately, additional assumptions are required for such an upper bound,
because (as earlier remarked) Lemma
\ref{lemma:convergence-for-almost-all-parameters} does not give {\em uniform}
convergence, merely universal convergence; with a large enough space of
hypotheses, the slowest pointwise convergence rates can be pushed arbitrarily
low.  For instance, let $\overline{x_1^t}$ be the distribution on
$\ObservableSpace^{\infty}$ which assigns probability 1 to endless repetitions
of $x_1^t$; clearly, under this distribution, seeing $X_1^t=x_1^t$ is almost
certain.  If such measures fall within the support of $\Posterior_0$, they will
dominate the likelihood, even though $\divrate(\overline{x_1^t}) = \infty$
under all but very special circumstances (e.g., $\TrueMeasure =
\overline{x_1^t}$).  Generically, then, the likelihood and the posterior weight
of $\overline{x_1^t}$ will rapidly plummet at times $T > t$.  To ensure
convergence of the posterior, overly-flexible measures like the family of
$\overline{x_1^t}$'s must be either excluded from the support of $\Posterior_0$
(possibly because they are excluded from $\ParameterSpace$), or be assigned so
little prior weight that they do not end up dominating the integrated
likelihood, or the posterior must concentrate on them.

\subsection{Convergence of Posterior Density via Control of the Integrated
  Likelihood}
\label{sec:control-of-prior}

The next two lemmas tell us that sets in $\ParameterSpace$ of
exponentially-small prior measure make vanishingly small contributions to the
integrated likelihood, and so to the posterior.  They do not require
assumptions beyond those used so far, but their application will.

\begin{lemma}
  Make Assumptions
  \ref{assumption:measurable-likelihood}--\ref{assumption:relative-AEP-holds},
  and chose a sequence of sets $B_t \subset \ParameterSpace$ such that, for all
  sufficiently large $t$, $\Posterior_0(B_t) \leq \alpha \myexp{-t \beta}$ for
  some $\alpha, \beta >0$.  Then, almost surely,
  \begin{equation}
    \prioravg{\LikelihoodRatio_t B_t} \leq \myexp{-t \beta/2}
  \label{eqn:vanishing-sets-make-vanishing-contributions}
  \end{equation}
  for all but finitely many $t$.
  \label{lemma:vanishing-sets-make-vanishing-contributions}
\end{lemma}
%lemma:4

\textsc{Proof:} By Markov's inequality.  First, use Fubini's theorem and the
chain rule for Radon-Nikodym derivatives to calculate the expectation value of
the ratio.
\begin{eqnarray*}
\Expect{\prioravg{\LikelihoodRatio_t B_t}}
& = &  \int_{\mathcal{X}^t}{d\TrueMeasure(x_1^t) \int_{B_t}{d\Posterior_0(\parametervalue) \LikelihoodRatio_t(\parametervalue) }}\\
& = & \int_{B_t}{d\Posterior_0(\parametervalue) \int_{\mathcal{X}^n}{d\TrueMeasure(x_1^t) \frac{d\ModelMeasure}{d\TrueMeasure}(x_1^t)}}\\
& = & \int_{B_t}{d\Posterior_0(\parametervalue) \int_{\mathcal{X}^t}{d\ModelMeasure(x_1^t)}}\\
& = & \int_{B_t}{d\Posterior_0(\parametervalue)}\\
& = & \Posterior_0(B_t)
\end{eqnarray*}
Now apply Markov's inequality:
\begin{eqnarray*}
\TrueMeasure\left\{x_1^t:~ \prioravg{\LikelihoodRatio_t B_t} > \myexp{-t \beta/2}\right\} & \leq & \myexp{t \beta/2} \Expect{\prioravg{\LikelihoodRatio_t B_t}}\\
& = & \myexp{t \beta/2} \Posterior_0(B_t)\\
& \leq & \alpha \myexp{-t \beta/2}
\end{eqnarray*}
for all sufficiently large $t$.  Since these probabilities are summable, the
Borel-Cantelli lemma implies that, with probability 1,
Eq.\ \ref{eqn:vanishing-sets-make-vanishing-contributions} holds for all but
finitely many $t$. $\Box$

The next lemma asserts a sequence of exponentially-small sets makes a
(logarithmically) negligible contribution to the posterior distribution,
provided the exponent is large enough compared to $\divrate(\ParameterSpace)$.

\begin{lemma}
  Let $B_t$ be as in the previous lemma.  If $\beta >
  2\divrate(\ParameterSpace)$, then
  \begin{equation}
    \frac{\prioravg{\LikelihoodRatio_t B_t^c}}{\prioravg{\LikelihoodRatio_t}} \rightarrow 1
  \end{equation}
  \label{lemma:integrated-likelihood-mostly-from-good-set}
\end{lemma}
%lemma:5

\textsc{Proof:} Begin with the likelihood integrated over $B_t$ rather than its
complement, and apply Lemmas \ref{lemma:integrated-likelihood-tracks-the-best}
and \ref{lemma:vanishing-sets-make-vanishing-contributions}: for any $\epsilon
> 0$
\begin{eqnarray}
  \frac{\prioravg{B_t \LikelihoodRatio_t}}{\prioravg{\LikelihoodRatio_n}}
  & \leq & \frac{\myexp{-t \beta/2}}{\myexp{-t[\divrate(\ParameterSpace) + \epsilon]}}\\
  & = & \myexp{t[\epsilon + \divrate(\ParameterSpace) - \beta/2]}
\label{eqn:integrated-likelihood-from-bad-set}
\end{eqnarray}
provided $t$ is sufficiently large.  If $\beta > 2 \divrate(\ParameterSpace)$,
this bound can be made to go to zero as $t\rightarrow\infty$ by taking
$\epsilon$ to be sufficiently small.  Since
\[
\prioravg{\LikelihoodRatio_t} = \prioravg{B_t^c \LikelihoodRatio_t} + \prioravg{B_t \LikelihoodRatio_t}
\]
it follows that
\[
\frac{\prioravg{B_t^c \LikelihoodRatio_t}}{\prioravg{\LikelihoodRatio_t}} \rightarrow 1
\]
$\Box$

\begin{lemma}
  Make Assumptions
  \ref{assumption:measurable-likelihood}--\ref{assumption:relative-AEP-holds},
  and take any set $G$ on which the convergence in Eq.\ \ref{eqn:relative-AEP}
  is uniform and where $\Posterior_0(G) > 0$.  Then, $\TrueMeasure$-a.s.,
  \begin{equation}
    \limsup_{t\rightarrow\infty}{\frac{1}{t}\log{\prioravg{G \LikelihoodRatio_t}}} \leq -\divrate(G)
  \end{equation}
  \label{lemma:lim-sup-of-integrated-likelihood-on-good-sets}
\end{lemma}
%lemma:6

\textsc{Proof:} Pick any $\epsilon > 0$.  By the hypothesis of uniform
convergence, there almost surely exists a $T(\epsilon)$ such that, for all $t
\geq T(\epsilon)$ and for all $\parametervalue \in G$,
$t^{-1}\log{\LikelihoodRatio_t(\parametervalue)} \leq
-\divrate(\parametervalue)+\epsilon$.  Hence
\begin{eqnarray}
t^{-1}\log{\prioravg{G \LikelihoodRatio_t}} & = & t^{-1}\log{\prioravg{G \myexp{\log{\LikelihoodRatio_t}}}}\\
& \leq & t^{-1}\log{\prioravg{G \myexp{t[-\divrate + \epsilon]}}}\\
& = & \epsilon + t^{-1}\log{\prioravg{G \myexp{-t\divrate}}} \label{eqn:integrated-likelihood-on-good-set-broken-up}
\end{eqnarray}
Let $\Posterior_{0|G}$ denote the probability measure formed by conditioning
$\Posterior_0$ to be in the set $G$.  Then
\[
\prioravg{G z} = \Posterior_0(G)\int_{G}{d\Posterior_{0|G}(\parametervalue) z(\parametervalue)}
\]
for any integrable function $z$.  Apply this to the last term from Eq.\
\ref{eqn:integrated-likelihood-on-good-set-broken-up}.
\[
\log{\prioravg{G \myexp{-t\divrate}}}
= \log{\Posterior_0(G)} + \log{\int_{G}{d\Posterior_{0|G}(\parametervalue)\myexp{-t\divrate(\parametervalue)}}}
\]
The second term on the right-hand side is the cumulant generating function of
$-\divrate(\parametervalue)$ with respect to $\Posterior_{0|G}$, which turns
out (cf.\ \citealp{Berk-consistency}) to have exactly the right behavior as $t \rightarrow\infty$.  
\begin{eqnarray}
\nonumber \frac{1}{t}\log{\int_{G}{d\Posterior_{0|G}(\parametervalue)\myexp{-t\divrate(\parametervalue)}}}  & = & \frac{1}{t}\log{\int_{G}{d\Posterior_{0|G}(\parametervalue){\left|\myexp{-\divrate(\parametervalue)}\right|}^t}}\\
\nonumber & = & \frac{1}{t}\log{{\left({\left(\int_{G}{d\Posterior_{0|G}(\parametervalue){\left|\myexp{-\divrate(\parametervalue)}\right|}^t} \right)}^{1/t} \right)}^t}\\
\nonumber & = & \frac{1}{t}\left(t \log{ {\left\|\myexp{-\divrate(\parametervalue)} \right\|}_{t,\Posterior_{0|G}}} \right)\\
& = & \log{ {\left\|\myexp{-\divrate(\parametervalue)} \right\|}_{t,\Posterior_{0|G}}}
\label{eqn:generating-function-divrate}
\end{eqnarray}
Since $\divrate(\parametervalue) \geq 0$, $\myexp{-\divrate(\parametervalue)}
\leq 1$, and the $L_p$ norm of the latter will grow towards its $L_{\infty}$
norm as $p$ grows.  Hence, for sufficiently large $t$,
\begin{eqnarray}
\nonumber \log{ {\left\|\myexp{-\divrate(\parametervalue)} \right\|}_{t,\Posterior_{0|G}}}
\nonumber & \leq & \log{ {\left\|\myexp{-\divrate(\parametervalue)} \right\|}_{\infty,\Posterior_{0|G}}} + \epsilon\\
\nonumber & = & -\essinf_{\parametervalue\in G}{\divrate(\parametervalue)} +\epsilon\\
& = & -\divrate(G) + \epsilon
\label{eqn:generating-function-vs-set-divrate}
\end{eqnarray}
where the next-to-last step uses the monotonicity of $\log$ and $\exp$.

Putting everything together, we have that, for any $\epsilon >0$ and all
sufficiently large $t$,
\begin{equation}
t^{-1}\log{\prioravg{G \LikelihoodRatio_t}} \leq -\divrate(G) + 2\epsilon +\frac{\log{\Posterior_0(G)}}{t}
\label{eqn:explicit-upper-bound-on-integrated-likelihood-for-large-t}
\end{equation}
Hence the limit superior of the left-hand side is at most $-\divrate(G)$.
$\Box$

\begin{lemma}
  Under Assumption
  \ref{assumption:measurable-likelihood}--\ref{assumption:sufficiently-rapid-convergence},
  \begin{equation}
    \limsup_{t\rightarrow\infty}{\frac{1}{t}\log{\prioravg{\LikelihoodRatio_t}}} \leq -\divrate(\ParameterSpace)
  \end{equation}
  \label{lemma:lim-sup-of-integrated-likelihood}
\end{lemma}
%lemma:7

\textsc{Proof:} By Lemma
\ref{lemma:integrated-likelihood-mostly-from-good-set},
\[
\lim_{t\rightarrow\infty}{\frac{\prioravg{\LikelihoodRatio_t}}{\prioravg{G_t \LikelihoodRatio_t}}} = 1
\]
implying that
\[
\lim_{t\rightarrow\infty}{\log{\prioravg{\LikelihoodRatio_t}} - \log{\prioravg{G_t \LikelihoodRatio_t}}} = 0
\]
so for every $\epsilon > 0$, for $t$ large enough
\[
\log{\prioravg{\LikelihoodRatio_t}} \leq \epsilon/3 + \log{\prioravg{G_t \LikelihoodRatio_t}} 
\]
Consequently, again for large enough $t$,
\[
\frac{1}{t}\log{\prioravg{\LikelihoodRatio_t}} \leq \epsilon/3t + \frac{1}{t} \log{\prioravg{G_t \LikelihoodRatio_t}} 
\]

Now, for each set $G$, for every $\epsilon > 0$, if $t \geq
\convtime(G,\epsilon/3)$ then
\[
\frac{1}{t} \log{\prioravg{G \LikelihoodRatio_t}} \leq -\divrate(G) + \epsilon/3
\]
by Lemma \ref{lemma:lim-sup-of-integrated-likelihood-on-good-sets}.  By
 Assumption
\ref{assumption:sufficiently-rapid-convergence}, $t \geq
\convtime(G_t,\epsilon/3)$ for all sufficiently large $t$.  Hence
\[
\frac{1}{t}\log{\prioravg{\LikelihoodRatio_t}} \leq -\divrate(G_t) +
\epsilon/3t + \epsilon/3
\]
for all $\epsilon > 0$ and all $t$ sufficiently large.  Since, by Assumption
\ref{assumption:exponentially-egorov}, $\divrate(G_t) \rightarrow
\divrate(\ParameterSpace)$, for every $\epsilon >0$, $\divrate(G_t)$ is within
$\epsilon/3$ of $\divrate(\ParameterSpace)$ for large enough $t$, so
\[
\frac{1}{t}\log{\prioravg{\LikelihoodRatio_t}} \leq -\divrate(\ParameterSpace) + \epsilon/3t + \epsilon/3 + \epsilon/3
\]
Thus, for every $\epsilon > 0$, then we have that
\[
\frac{1}{t}\log{\prioravg{\LikelihoodRatio_t}} \leq -\divrate(\ParameterSpace)
+ \epsilon
\]
for large enough $t$, or, in short,
\[
\limsup_{t\rightarrow\infty}{\frac{1}{t}\log{\prioravg{\LikelihoodRatio_t}}} \leq -\divrate(\ParameterSpace)
\]
$\Box$

\begin{lemma}
  Under Assumptions
  \ref{assumption:measurable-likelihood}--\ref{assumption:sufficiently-rapid-convergence},
  if $\Posterior_0(\InfinitySet) = 0$, then
  \begin{equation}
    \frac{1}{t}\log{\prioravg{\LikelihoodRatio_t}} \rightarrow -\divrate(\ParameterSpace)
  \end{equation}
  almost surely.
  \label{lemma:time-avg-of-log-integrated-likelihood}
\end{lemma}
%lemma:8

\textsc{Proof}: Combining Lemmas
\ref{lemma:integrated-likelihood-tracks-the-best} and
\ref{lemma:lim-sup-of-integrated-likelihood},
\[
-\divrate(\ParameterSpace)
\leq \liminf_{t\rightarrow\infty}{\frac{1}{t}\log{\prioravg{\LikelihoodRatio_t}}}
\leq \limsup_{t\rightarrow\infty}{\frac{1}{t}\log{\prioravg{\LikelihoodRatio_t}}}
\leq -\divrate(\ParameterSpace)
\]
$\Box$

The standard version of Egorov's theorem concerns sequences of finite
measurable functions converging pointwise to a finite measurable limiting
function.  However, the proof is easily adapted to an infinite limiting
function.

\begin{lemma}
  Let $f_t(\parametervalue)$ be a sequence of finite, measurable functions,
  converging to $\infty$ almost everywhere ($\Posterior_0$) on $I$.  Then for
  each $\epsilon > 0$, there exists a possibly-empty $B \subset I$ such that
  $\Posterior_0(B) < \epsilon$, and the convergence is uniform on $I \setminus
  B$.
  \label{lemma:egorov-with-infinite-limit}
\end{lemma}
%lemma:9

\textsc{Proof:} Parallel to the usual proof of Egorov's theorem.  Begin by
removing the measure-zero set of points on which pointwise convergence fails;
for simplicity, keep the name $I$ for the remaining set.  For each natural
number $t$ and $k$, let $B_{t,k} \equiv \left\{ \parametervalue \in I:~
  f_t(\parametervalue) < k\right\}$ --- the points where the function fails to
be at least $k$ by step $t$.  Since the limit of $f_t$ is $\infty$ everywhere
on $I$, each $\parametervalue$ has a last $t$ such that $f_t(\parametervalue) <
k$, no matter how big $k$ is.  Hence $\bigcap_{t=1}^{\infty}{B_{t,k}} =
\emptyset$.  By continuity of measure, for any $\delta > 0$, there exists an
$n$ such that $\Posterior_0(B_{t,k}) < \delta$ if $t \geq n$.  Fix $\epsilon$
as in the statement of the lemma, and set $\delta = \epsilon 2^{-k}$.  Finally,
set $B = \bigcup_{k=1}^{\infty}{B_{n,k}}$.  By the union bound,
$\Posterior_0(B) \leq \epsilon$, and by construction, the rate of convergence
to $\infty$ is uniform on $I\setminus B$. $\Box$

\begin{lemma}
  The conclusion of Lemma \ref{lemma:time-avg-of-log-integrated-likelihood} is
  unchanged if $\Posterior_0(\InfinitySet) > 0$.
  \label{lemma:time-avg-of-log-integrated-likelihood-with-infinity}
\end{lemma}
%lemma:10

\textsc{Proof:} The integrated likelihood ratio can be divided into two parts,
one from integrating over $\InfinitySet$ and one from integrating over its
complement.  Previous lemmas have established that the latter is upper bounded,
in the long run, by a quantity which is
$O(\myexp{-\divrate(\ParameterSpace)t})$.  We can use Lemma
\ref{lemma:egorov-with-infinite-limit} to divide $\InfinitySet$ into a sequence
of sub-sets, on which the convergence is uniform, and hence on which the
integrated likelihood shrinks faster than any exponential function, and
remainder sets, of prior measure no more than $\alpha \myexp{-n \beta}$, on
which the convergence is less than uniform (i.e., slow).  If we ensure that
$\beta > 2 \divrate(\ParameterSpace)$, however, by Lemma
\ref{lemma:integrated-likelihood-mostly-from-good-set} the remainder sets'
contributions to the integrated likelihood is negligible in comparison to that
of $\ParameterSpace \setminus \InfinitySet$.  Said another way, if there are
alternatives which a consistent test would rule out at a merely exponential
rate, those which would be rejected at a supra-exponential rate end up making
vanishingly small contributions to the integrated likelihood. $\Box$

\begin{theorem}
  Under Assumptions
  \ref{assumption:measurable-likelihood}--\ref{assumption:sufficiently-rapid-convergence},
  for all $\parametervalue \in \ParameterSpace$ where
  $\posteriordensity_0(\parametervalue) > 0$,
  \begin{equation}
    \lim_{t\rightarrow\infty}{\frac{1}{t}\log{\posteriordensity_t(\parametervalue)}} = -\Action(\parametervalue)
  \end{equation}
  with probability 1.
  \label{theorem:long-run-fitness}
\end{theorem}
%theorem:2

\textsc{Proof:} Theorem \ref{thm:upper-bound-on-long-run-fitness} says that,
for all $\parametervalue$,
\[
\limsup_{t\rightarrow\infty}{\frac{1}{t}\log{\posteriordensity_t(\parametervalue)}} \leq -\Action(\parametervalue)
\]
a.s., so there just needs to be a matching $\liminf$.  Pick any $\epsilon > 0$.
By Assumption \ref{assumption:relative-AEP-holds}, it's almost certain that,
for all sufficiently large $t$,
\[
\frac{1}{t}\log{\LikelihoodRatio_t(\parametervalue)} \geq -\divrate(\theta) - \epsilon/2
\]
while by Lemma \ref{lemma:time-avg-of-log-integrated-likelihood-with-infinity},
it's almost certain that for all sufficiently large $t$,
\[
\frac{1}{t}\log{\prioravg{\LikelihoodRatio_t}} \leq -\divrate(\ParameterSpace) + \epsilon/2
\]
Combining these as in the proof of Theorem
\ref{thm:upper-bound-on-long-run-fitness}, it's almost certain that for all
sufficiently large $t$
\[
\frac{1}{t}\log{\posteriordensity_t(\parametervalue)} \geq \divrate(\ParameterSpace) - \divrate(\parametervalue) - \epsilon
\]
so
\[
\liminf_{t\rightarrow\infty}{\frac{1}{t}\log{\posteriordensity_t(\parametervalue)}} \geq \divrate(\ParameterSpace) - \divrate(\parametervalue) = -\Action(\parametervalue)
\]
$\Box$

\subsection{Convergence and Large Deviations of the Posterior Measure}
\label{sec:convergence-and-ldp}

Adding Assumption \ref{assumption:good-sets-are-good-even-for-subsets} to those
before it implies that the posterior measure concentrates on sets $A \subset
\ParameterSpace$ where $\divrate(A) = \divrate(\ParameterSpace)$.

\begin{theorem}
  Make Assumptions
  \ref{assumption:measurable-likelihood}--\ref{assumption:good-sets-are-good-even-for-subsets}.
  Pick any set $A\in\ParameterSigmaField$ where $\Posterior_0(A) > 0$ and
  $\divrate(A) > \divrate(\ParameterSpace)$.  Then $\Posterior_t(A) \rightarrow
  0$ a.s.
  \label{thm:convergence-of-measure}
\end{theorem}
%theorem:3

\textsc{Proof}: 
\begin{eqnarray*}
\Posterior_t(A) & = & \Posterior_t(A \cap G_t) + \Posterior_t(A \cap G^c_t)\\
& \leq & \Posterior_t(A \cap G_t) + \Posterior_t(G^c_t)
\end{eqnarray*}
The last term is easy to bound.  From Eq.\
\ref{eqn:integrated-likelihood-from-bad-set} in the proof of Lemma
\ref{lemma:integrated-likelihood-mostly-from-good-set},
\begin{eqnarray}
\nonumber \Posterior_t(G^c_t) & = & \frac{\prioravg{\LikelihoodRatio_t G_t^c}}{\prioravg{\LikelihoodRatio_t}}\\
& \leq & \myexp{t[\epsilon + \divrate(\ParameterSpace) - \beta/2]}
\label{eqn:posterior-of-remainder}
\end{eqnarray}
for any $\epsilon > 0$, for all sufficiently large $t$, almost surely.  Since
$\beta > 2\divrate(\ParameterSpace)$, the whole expression $\rightarrow 0$ as
$t \rightarrow \infty$.

To bound $\Posterior_t(A \cap G_t)$, reasoning as in the proof of Lemma
\ref{lemma:lim-sup-of-integrated-likelihood}, but invoking Assumption
\ref{assumption:good-sets-are-good-even-for-subsets}, leads to the conclusion
that, for any $\epsilon > 0$, with probability 1,
\[
\frac{1}{t}\log{\prioravg{\LikelihoodRatio_t (A \cap G_t)}} \leq -\divrate(A) +\epsilon
\]
for all sufficiently large $t$.  Recall that by Lemma
\ref{lemma:integrated-likelihood-tracks-the-best}, for all $\epsilon > 0$ it's
almost sure that
\[
\frac{1}{t}\log{\prioravg{\LikelihoodRatio_t}} \geq -\divrate(\ParameterSpace) - \epsilon
\]
for all sufficiently large $n$.  Hence for every $\epsilon > 0$, it's almost
certain that for all sufficiently large $t$,
\begin{equation}
\Posterior_t(A \cap G_t) \leq \myexp{t[\divrate(\ParameterSpace) - \divrate(A) + 2\epsilon]}
\label{eqn:posterior-of-A-and-good-set}
\end{equation}
Since $\divrate(A) > \divrate(\ParameterSpace)$, by picking $\epsilon$ small
enough the right hand side goes to zero.  $\Box$

The proof of the theorem provides an exponential {\em upper} bound on the
posterior measure of sets where $\divrate(A) > \divrate(\ParameterSpace)$.  In
fact, even without the final assumption needed for the theorem, there is an
exponential {\em lower} bound on that posterior measure.

\begin{lemma}
  Make Assumption
  \ref{assumption:measurable-likelihood}--\ref{assumption:sufficiently-rapid-convergence},
  and pick a set $A \in \ParameterSigmaField$ with $\Posterior_0(A) > 0$.  Then
  \begin{equation}
    \liminf_{t\rightarrow\infty}{\frac{1}{t}\log{\Posterior_t(A)}} \geq -\Action(A)
  \end{equation}
\end{lemma}
%lemma:11

\textsc{Proof:} Reasoning as in the proof of Lemma
\ref{lemma:integrated-likelihood-tracks-the-best}, it is easy to see that
\[
\liminf_{t\rightarrow\infty}{\frac{1}{t}\log{\prioravg{\LikelihoodRatio_t A}}} \geq -\divrate(A)
\]
and by Lemma \ref{lemma:lim-sup-of-integrated-likelihood},
\[
\limsup_{t\rightarrow\infty}{\frac{1}{t} \log{\prioravg{\LikelihoodRatio_t}}} \leq -\divrate(\ParameterSpace) 
\]
hence
\begin{eqnarray*}
\liminf_{t\rightarrow\infty}{\frac{1}{t}\log{\Posterior_t(A)}} & = & \liminf_{t\rightarrow\infty}{\frac{1}{t}\log{\frac{\prioravg{\LikelihoodRatio_t A}}{\prioravg{\LikelihoodRatio_t}}}} \\
& \geq &  - \divrate(A) + \divrate(\ParameterSpace)
\end{eqnarray*}
$\Box$

\begin{theorem}
  Under the conditions of Theorem \ref{thm:convergence-of-measure}, if $A\in\ParameterSigmaField$
  is such that
  \begin{equation}
  -\limsup{t^{-1}\log{\Posterior_0(A \cap G_t^c)}} = \beta^{\prime} \geq
  2\divrate(A)
  \end{equation}
  then
  \begin{equation}
    \lim_{t\rightarrow\infty}{\frac{1}{t}\log{\Posterior_t(A)}} = \divrate(\ParameterSpace) - \divrate(A)
    \label{eqn:ldp}
  \end{equation}
  In particular, this holds whenever $2\divrate(A) < \beta$ or $A \subset
  \bigcap_{k=n}^{\infty}{G_k}$ for some $n$.
  \label{thm:ldp}
\end{theorem}
%theorem:4

\textsc{Proof:} Trivially,
\[
\frac{1}{t}\log{\Posterior_t(A)} = \frac{1}{t}\log{\Posterior_t(A \cap G_t) + \Posterior(A \cap G_t^c)}
\]
From Eq.\ \ref{eqn:posterior-of-A-and-good-set} from the proof of Theorem
\ref{thm:convergence-of-measure}, we know that, for any $\epsilon > 0$,
\[
\Posterior_t(A \cap G_t) \leq \myexp{t[\divrate(\ParameterSpace) - \divrate(A) + \epsilon]}
\]
a.s. for sufficiently large $t$.  
On the other hand, under the hypothesis of the theorem, the proof of Eq.\
\ref{eqn:posterior-of-remainder} can be imitated for $\Posterior_t(A \cap
G_t^c)$, with the conclusion that, for all $\epsilon > 0$,
\[
\Posterior_t(A \cap G_t^c) \leq \myexp{t[\divrate(\ParameterSpace) - \beta^{\prime}/2 + \epsilon]}
\]
again a.s. for sufficiently large $t$.  Since $\beta^{\prime}/2 > \divrate(A)$,
Eq.\ \ref{eqn:ldp} follows.

Finally, to see that this holds for any $A$ where $\divrate(A) < \beta/2$,
observe that we can always upper bound $\Posterior_t(A \cap G_t^c)$ by
$\Posterior_t(G_t^c)$, and the latter goes to zero with rate at least
$-\beta/2$. $\Box$

{\em Remarks:} Because $\divrate(A)$ is the essential infimum of
$\divrate(\theta)$ on the set $A$, as the set shrinks $\divrate(A)$ grows.
Sets where $\divrate(A)$ is much larger than $\divrate(\ParameterSpace)$ tend
accordingly to be small.  The difficulty is that the sets $G_t^c$ are also
small, and conceivably overlap so heavily with $A$ that the integral of the
likelihood over $A$ is dominated by the part coming from $A \cap G_t^c$.
Eventually this will shrink towards zero exponentially, but perhaps only at the
comparatively slow rate $\divrate(\ParameterSpace) - \beta/2$, rather than the
faster rate $\divrate(\ParameterSpace) - \divrate(A)$ attained on the
well-behaved part $A \cap G_t$.

Theorem \ref{thm:ldp} is close to, but not quite, a large deviations principle
on $\ParameterSpace$.  We have shown that the posterior probability of any
arbitrary set $A$ where $\Action(A) > 0$ goes to zero with an exponential rate
at least equal to
\begin{equation}
\beta/2 \wedge \essinf_{\parametervalue \in A}{\Action(\parametervalue)} = \essinf_{\parametervalue \in A}{\beta/2 \wedge \Action(\parametervalue)}
\end{equation}
But in a true LDP, the rate would have to be an infimum, not just an essential
infimum, of a point-wise rate function.  This deficiency could be removed by
means of additional assumptions on $\Posterior_0$ and
$\divrate(\parametervalue)$.

Ref.\ \citep{Eichelsbacher-Ganesh-moderate-dev-for-posteriors} obtains proper
large and even moderate deviations principles, but for the location of
$\Posterior_t$ in the space $\mathcal{M}_1(\ParameterSpace)$ of all
distributions on $\ParameterSpace$, rather than on $\ParameterSpace$ itself.
Essentially, they use the assumption of IID sampling, which makes the posterior
a function of the empirical distribution, to leverage the LDP for the latter
into an LDP for the former.  This strategy may be more widely applicable but
goes beyond the scope of this paper.  Papangelou
\citep{Papangelou-Bayesian-estimation-of-higher-order-markov}, assuming that
$\ParameterSpace$ consists of discrete-valued Markov chains of arbitrary order
and $\TrueMeasure$ is in the support of the prior, and using methods similar to
those in Appendix \ref{appendix:verification}, derives a result which is
closely related to Theorem \ref{thm:ldp}.  In fact, fixing the sets $G_t$ as in
Appendix \ref{appendix:verification}, Theorem \ref{thm:ldp} implies the theorem
of \citep{Papangelou-Bayesian-estimation-of-higher-order-markov}.

\subsection{Generalization Performance}
\label{sec:generalization-performance}

Lemma \ref{lemma:time-avg-of-log-integrated-likelihood-with-infinity} shows
that, in hindsight, the Bayesian learner does a good job of matching the data:
the log integrated likelihood ratio per time-step approaches
$-\divrate(\ParameterSpace)$, the limit of values attainable by individual
hypotheses within the support of the prior.  This leaves open, however, the
question of the prospective or generalization performance.

What we want is for the posterior predictive distribution
$\PredictiveMeasure^t$ to approach the true conditional distribution of future
events $\TrueMeasure^t$, but we cannot in general hope for the convergence to
be complete, since our models are mis-specified.  The next theorem uses
$\divrate(\ParameterSpace)$ to put an upper bound on how far the posterior
predictive distribution can remain from the true predictive distribution.

\begin{theorem}
  Under Assumptions
  \ref{assumption:measurable-likelihood}--\ref{assumption:good-sets-are-good-even-for-subsets}, with probability 1,
  \begin{eqnarray}
    \label{eqn:hellinger-bound} \limsup_{t\rightarrow\infty}{\rho^2_{H}(\TrueMeasure^t,\PredictiveMeasure^t)} & \leq & \divrate(\ParameterSpace)\\
    \limsup_{t\rightarrow\infty}{\rho_{TV}^2(\TrueMeasure^t,\PredictiveMeasure^t)} & \leq & 4\divrate(\ParameterSpace)
    \label{eqn:tv-bound}
  \end{eqnarray}
  where $\rho_{H}$ and $\rho_{TV}$ are, respectively, the Hellinger and total
  variation metrics.
\label{thm:generalization-performance}
\end{theorem}
%theorem:5

\textsc{Proof:} Recall the well-known inequalities relating Hellinger distance
to Kullback-Leibler divergence on the one side and to total variation
distance on the other \citep{Ghosh-Ramamoorthi}: for any two distributions $P$
and $Q$,
\begin{eqnarray}
\label{eqn:hellinger-entropy-inequality} \rho_{H}^2(P,Q) & \leq & D(P\|Q)\\
\rho_{TV}(P,Q) & \leq & 2\rho_{H}(P,Q)
\label{eqn:hellinger-tv-inequality}
\end{eqnarray}
It's enough to prove Eq.\ \ref{eqn:hellinger-bound}, and Eq.\
\ref{eqn:tv-bound} then follows from Eq.\ \ref{eqn:hellinger-tv-inequality}.

Abbreviate $\rho_{H}(\TrueMeasure^t,\ModelMeasure^t)$ by
$\rho_H(t,\parametervalue)$.  Pick any $\epsilon > \divrate(\ParameterSpace)$,
and say that $A_\epsilon = \left\{\parametervalue: \rho^2_H(t,\parametervalue)
  > \epsilon\right\}$.  By convexity and Jensen's inequality,
\begin{eqnarray*}
\rho^2_{H}(\TrueMeasure^t,\PredictiveMeasure^t)
& \leq & \int_{\ParameterSpace}{\rho^2_H(t,\parametervalue) d\Posterior_n(\parametervalue)}\\
& = & \int_{A_\epsilon^c}{\rho^2_H(t,\parametervalue) d\Posterior_n(\parametervalue)} + \int_{A_\epsilon}{\rho^2_H(t,\parametervalue) d\Posterior_n(\parametervalue)}\\
& = & \epsilon \Posterior_t(A_{\epsilon}^c) + \sqrt{2}\Posterior_t(A_{\epsilon})
\end{eqnarray*}
By Eq.\ \ref{eqn:hellinger-entropy-inequality}, $d(\parametervalue) >
\rho^2_H(t,\parametervalue)$. Thus $\divrate(A_{\epsilon}) \geq \epsilon$, and
$\epsilon > \divrate(\parametervalue)$ so, by Theorem
\ref{thm:convergence-of-measure}, $\Posterior_t(A_{\epsilon}) \rightarrow 0$
a.s.  Hence
\[
\rho^2_{H}(\TrueMeasure^t,\PredictiveMeasure^t) \leq \epsilon
\]
eventually almost surely.  Since this holds for any $\epsilon >
\divrate(\ParameterSpace)$, Eq. \ref{eqn:hellinger-bound} follows.  $\Box$

{\em Remark}: It seems like it should be possible to prove a similar result for
the divergence rate of the predictive distribution, namely that
\[
\limsup_{t\rightarrow\infty}{\divrate(\Posterior_t)} \leq \divrate(\ParameterSpace)
\]
but it would take a different approach, because $\divrate(\cdot)$ has no upper
bound, and the posterior weight of the high-divergence regions might decay too
slowly to compensate for this.

\subsection{Rate of Convergence}
\label{sec:rate-of-convergence}

Recall that $N_{\epsilon}$ was defined as the set of all $\parametervalue$ such
that $\divrate(\parametervalue) \leq \divrate(\ParameterSpace) + \epsilon$.
(This is measurable by Assumption \ref{assumption:divergence-rates-exist}.)
The set $N_{\epsilon}^c$ thus consists of all hypotheses whose divergence rate
is more than $\epsilon$ above the essential infimum
$\divrate(\ParameterSpace)$.  For any $\epsilon > 0$,
$\Posterior_t(N_{\epsilon}^c) \rightarrow 0$ a.s., by Theorem
\ref{thm:convergence-of-measure}, and for sufficiently small $\epsilon$,
$\lim_{t\rightarrow\infty}{t^{-1} \log{\Posterior_t(N_{\epsilon}^c)}} =
-\epsilon$ a.s., by Theorem \ref{thm:ldp}.  For such sets, in other words, for
any $\delta > 0$, it's almost certain that for all sufficiently large $t$,
\begin{equation}
  \Posterior_t(N_{\epsilon}^c) \leq \myexp{-t(\epsilon -\delta)}
  \label{eqn:exponental-posterior-shrinkage-of-suboptimal-sets}
\end{equation}
Now consider a non-increasing positive sequence $\epsilon_t \rightarrow 0$.
Presumably if $\epsilon_t$ decays slowly enough,
$\Posterior_t(N^c_{\epsilon_t})$ will still go to zero, even though the sets
$N^c_{\epsilon_t}$ are non-decreasing.  Examination of Eq.\
\ref{eqn:exponental-posterior-shrinkage-of-suboptimal-sets} suggests, naively,
that this will work if $t\epsilon_t \rightarrow \infty$, i.e., if the decay of
$\epsilon_t$ is strictly sublinear.  This is correct under an additional
condition, similar to Assumption
\ref{assumption:sufficiently-rapid-convergence}.

\begin{theorem}
  Make assumptions
  \ref{assumption:measurable-likelihood}--\ref{assumption:good-sets-are-good-even-for-subsets},
  and pick a positive sequence $\epsilon_t$ where $\epsilon_t \rightarrow 0$,
  $t\epsilon_t \rightarrow\infty$.  If, for each $\delta > 0$, 
  \begin{equation}
    \tau(G_t \cap N^c_{\epsilon_t},\delta) \leq t
    \label{eqn:sufficiently-rapid-convergence-for-rates}
  \end{equation}
  eventually almost surely, then
  \begin{equation}
    \Posterior_t(N_{\epsilon_t}) \rightarrow 1
  \end{equation}
  with probability 1.
  \label{thm:rate-of-convergence}
\end{theorem}

\textsc{Proof}: By showing that $\Posterior_t(N^c_{\epsilon_t}) \rightarrow 0$
a.s.  Begin by splitting the sets into the parts inside the $G_t$, say 
$U_t$, and the parts outside:
\begin{eqnarray*}
\Posterior_t(N^c_{\epsilon_t}) &= & \Posterior_t(N^c_{\epsilon_t} \cap G_t) + \Posterior_t(N^c_{\epsilon_t} \cap G^c_t)\\
& \leq & \Posterior_t(U_t) + \Posterior_t(G^c_t)
\end{eqnarray*}
From Lemma \ref{lemma:vanishing-sets-make-vanishing-contributions}, the second
term $\rightarrow 0$ with probability 1, so for any $\eta_1 > 0$, it is $\leq
\eta_1$ eventually a.s.

Turning to the other term, Theorem \ref{thm:ldp} applies to $U_k$ for any
fixed $k$, so
\[
\lim_{t\rightarrow\infty}{t^{-1}\log{\Posterior_t(U_k)}} = \divrate(\ParameterSpace) - \divrate(U_k)
\]
(a.s.), implying, with Lemma
\ref{lemma:time-avg-of-log-integrated-likelihood-with-infinity}, that
\[
\lim_{t\rightarrow\infty}{t^{-1}\log{\prioravg{U_k \LikelihoodRatio_t}}} = -\divrate(U_k)
\]
(a.s.).  By Eq.\ \ref{eqn:sufficiently-rapid-convergence-for-rates}, for any
$\eta_2 > 0$,
\[
t^{-1}\log{\prioravg{U_t \LikelihoodRatio_t}} \leq -\divrate(U_t) + \eta_2
\]
eventually almost surely.  By Lemma \ref{lemma:time-avg-of-log-integrated-likelihood-with-infinity} and Bayes's rule, then,
\[
t^{-1}\log{\Posterior_t(U_t)} \leq \divrate(\ParameterSpace) - \divrate(U_t) + \eta_3
\]
eventually a.s., for any $\eta_3 > 0$.  Putting things back together, eventually
a.s.,
\begin{eqnarray*}
\Posterior_t(N^c_{\epsilon_t}) & \leq & \myexp{t(\divrate(\ParameterSpace) - \divrate(U_t) + \eta_2)} + \eta_1\\
& \leq & \myexp{t(-\epsilon_t + \eta_3)} + \eta_1
\end{eqnarray*}
Since $t\epsilon_t \rightarrow \infty$, the first term goes to zero, and since
$\eta_1$ can be as small as desired,
\[
\Posterior_t(N^c_{\epsilon_t}) \rightarrow 0
\]
almost surely.  $\Box$

The theorem lets us attain rates of convergence just slower than $t^{-1}$ (so
that $t\epsilon_t \rightarrow \infty$).  This matches existing results on rates
of posterior convergence for mis-specified models with IID data in
\cite[Corollary 5.2]{Zhang-from-epsilon-to-KL} ($t^{-1}$ in the Renyi
divergence) and in \citep{Kleijn-van-der-Vaart} ($t^{-1/2}$ in the Hellinger
distance; recall Eq.\ \ref{eqn:hellinger-entropy-inequality}), and for
correctly-specified non-IID models in
\citep{Ghosal-van-der-Vaart-posterior-convergence-non-iid} ($t^{-\alpha}$ for
suitable $\alpha < 1/2$, again in the Hellinger distance).

\subsection{Application of the Results to the Example}
\label{sec:results-on-example}

Because $\divrate(\ParameterSpace) = 0$, while $\divrate(\parametervalue) > 0$
everywhere, the behavior of the posterior is somewhat peculiar.  Every compact
set $K \subset \ParameterSpace$ has $\Action(K) > 0$, so by Theorem
\ref{thm:convergence-of-measure}, $\Posterior_t(K) \rightarrow 0$.  On the
other hand, $\Posterior_t(G_t) \rightarrow 1$ --- the sequence of good sets
contains models of increasingly high order, with increasingly weak constraints
on the transition probabilities, and this lets its posterior weight grow, even
though every individual compact set within it ultimately loses all weight.

In fact, each $G_t$ is a convex set, and $\divrate(\cdot)$ is a convex
function, so there is a unique minimizer of the divergence rate within each
good set.  Conditional on being within $G_t$, the posterior probability becomes
increasingly concentrated on neighborhoods of this minimizer, but the minimizer
itself keeps moving, since it can always be improved upon by increasing the
order of the chain and reducing some transition probabilities.  (Recall that
$\TrueMeasure$ gives probability 0 to sequences $010$, $01110$, etc., where the
block of 1's is of odd length, but $\ParameterSpace$ contains only chains with
strictly positive transition probabilities.)

Outside of the good sets, the likelihood is peaked around hypotheses which
provide stationary and smooth approximations to the $\overline{x_1^t}$
distribution that endlessly repeats the observed sequence to date. The
divergence rates of these hypotheses are however extremely high, so none of
them retains its high likelihood for very long.  ($\overline{x_1^t}$ is a
Markov chain of order $t$, but it is not in $\ParameterSpace$, since it's
neither stationary nor does it have strictly positive transition probabilities.
It can be made stationary, however, by assigning equal probability to each of
its $t$ states; this gives the data likelihood $1/t$ rather than 1, but that
still is vastly larger than the $O(-ct)$ log-likelihoods of better models.
(Recall that even the log-likelihood of the true distribution is only
$O(-\frac{2}{3}t)$.)  Allowing each of the $t$ states to have a probability $0
< \iota \ll 1$ of not proceeding to the next state in the periodic sequence is
easy and leads to only an $O(\iota t)$ reduction in the likelihood up to time
$t$.  In the long run, however, it means that the log-likelihood will be
$O(t\log{\iota})$.)  In any case, the total posterior probability of $G_t^c$ is
going to zero exponentially.

Despite --- or rather, because of ---  the fact that no point in
$\ParameterSpace$ is the {\em ne plus ultra} around which the posterior
concentrates, the conditions of Theorem \ref{thm:generalization-performance}
are met, and since $\divrate(\ParameterSpace) = 0$, the posterior predictive
distribution converges to the true predictive distribution in the Hellinger and
total variation metrics.  That is, the weird gyrations of the posterior do not
prevent us from attaining {\em predictive} consistency.  This is so even though
the posterior always gives the wrong answer to such basic questions as ``Is
$\TrueMeasure(X_{t}^{t+2} = 010) > 0$?'' --- inferences which in this case can
be made correctly through non-Bayesian methods
\citep{Ornstein-Weiss-how-sampling-reveals-a-process,CSSR-for-UAI}.

\section{Discussion}
\label{sec:conclusion}

The crucial assumptions were \ref{assumption:relative-AEP-holds},
\ref{assumption:exponentially-egorov} and
\ref{assumption:sufficiently-rapid-convergence}.  Together, these amount to
assuming that the time-averaged log likelihood ratio converges universally; to
fashioning a sieve, successively embracing regions of $\ParameterSpace$ where
the convergence is increasingly ill-behaved; and the hope that the prior weight
of the remaining bad sets can be bounded exponentially.

Using asymptotic equipartition in place of the law of large numbers is fairly
straightforward.  Both results belong to the general family of ergodic
theorems, which allow us to take sufficiently long sample paths as
representative of entire processes.  The unique a.s. limit in Eq.\
\ref{eqn:relative-AEP} can be replaced with a.s. convergence to a distinct
limit in each ergodic component of $\TrueMeasure$.  However, the notation gets
ugly, so the reader should regard $\divrate(\parametervalue)$ as that random
limit, and treat all subsequent results as relative to the ergodic
decomposition of $\TrueMeasure$.
(Cf. \citep{Gray-ergodic-properties,Debowski-ergodic-decomposition}.)  It may be
possible to weaken this assumption yet further, but it is hard to see how
Bayesian updating can succeed if the past performance of the likelihood is not
a guide to future results.

A bigger departure from the usual approach to posterior convergence may be
allowing $\divrate(\ParameterSpace) > 0$; this rules out posterior consistency,
to begin with.  More subtly, it requires $\beta > 2\divrate(\ParameterSpace)$.
This means that a prior distribution which satisfies the assumptions for one
value of $\TrueMeasure$ may not satisfy them for another, depending, naturally
enough, on just how mis-specified the hypotheses are, and how much weight the
prior puts on very bad hypotheses.  On the other hand, when
$\divrate(\ParameterSpace) = 0$, Theorem \ref{thm:generalization-performance}
implies predictive consistency, as in the example.

Assumption \ref{assumption:sufficiently-rapid-convergence} is frankly annoying.
It ensures that the log likelihood ratio converges fast enough, at least on
the good sets, that we can be confident that integrated likelihood of $G_t$ has
converged well by the time we want $G_t$ to start dominating the prior.  It was
shaped, however, to fill a hole in the proof of Lemma
\ref{lemma:lim-sup-of-integrated-likelihood} rather than more natural
considerations.  The result is that verifying the assumption in its present
form means proving the sub-linear growth rate of sequences of random last entry
times, and these times are not generally convenient to work with.  (Cf.\
Appendix \ref{appendix:verification}.)  It would be nice to replace it with a
bracketing or metric entropy condition, as in
\citep{Barron-Schervish-Wasserman,Zhang-from-epsilon-to-KL} or similar forms of
capacity control, as used in
\citep{Meir-nonparametric-time-series,Vidyasagar-on-learning-and-generalization}.
Alternately, the uniformly consistent test conditions widely employed in
Bayesian nonparametrics
\citep{Ghosh-Ramamoorthi,Xing-Ranneby-necessary-and-sufficient-exponential-consistency}
have been adapted the mis-specified setting by \citep{Kleijn-van-der-Vaart},
where the tests become reminiscent of the ``model selection tests'' used in
econometrics \citep{Vuong-testing-non-nested-hypotheses}.  Since the latter can
work for dynamical models \citep{Rivers-Vuong-model-selection-tests}, this
approach may also work here.  In any event, replacing Assumption
\ref{assumption:sufficiently-rapid-convergence} with more primitive,
comprehensible and easily-verified conditions seems a promising direction for
future work.

These results go some way toward providing a frequentist explanation of the
success of Bayesian methods in many practical problems.  Under these
conditions, the posterior is increasingly weighted towards the parts of
$\ParameterSpace$ which are closest (in the Kullback-Leibler sense) to the
data-generating process $\TrueMeasure$.  For a $\Posterior_t(A)$ to
persistently be much more or much less than $\approx \myexp{-t\Action(A)}$,
$\LikelihoodRatio(\parametervalue)$ must be persistently far from
$\myexp{-t\divrate(\parametervalue)}$, not just for isolated $\parametervalue
\in A$, but a whole positive-measure subset of them.  With a reasonably smooth
prior, this requires a run of bad luck amounting almost to a conspiracy.  From
this point of view, Bayesian inference amounts to introducing bias so as to
reduce variance, and then relaxing the bias.  Experience with frequentist
non-parametric methods shows this can work if the bias is relaxed sufficiently
slowly, which is basically what the assumptions here do.  As the example shows,
this can succeed as a predictive tactic without supporting substantive
inferences about the data-generating process.  However,
\ref{assumption:best-divergence-is-finite}--\ref{assumption:good-sets-are-good-even-for-subsets}
involve both the prior and the data-generating process, and so cannot be
verified using the prior alone.  For empirical applications, it would be nice
to have ways of checking them using sample data.

When $\divrate(\ParameterSpace) > 0$ and all the models are more or less wrong,
there is an additional advantage to averaging the models, as is done in the
predictive distribution.  (I owe the argument which follows to Scott Page;
cf. \citep{scotte-Difference}.)  With a convex loss function $\ell$, such as
squared error, Kullback-Leibler divergence, Hellinger distance, etc., the loss
of the predictive distribution $\ell(\Posterior_t)$ will be no larger than the
posterior-mean loss of the individual models
$\popavgt{\ell(\parametervalue)}{t}$.  For squared error loss, the difference
is equal to the variance of the models' predictions
\citep{Krogh-Vedelsby-neural-network-ensembles}.  For divergence, some algebra
shows that
\begin{equation}
\divrate(\Posterior_t) = \popavgt{\divrate(\parametervalue)}{t} + \popavgt{\Expect{\log{\frac{d\ModelMeasure}{d\PredictiveMeasure}}}}{t}
\end{equation}
where the second term on the RHS is again an indication of the diversity of the
models; the more different their predictions are, on the kind of data generated
by $\TrueMeasure$, the smaller the error of made by the mixture.  Having a
diversity of wrong answers can be as important as reducing the average error
itself.  The way to accomplish this is to give more weight to models which make
mostly good predictions, but make different mistakes.  This suggests that there
may actually be predictive benefits to having the posterior concentrate on a
set containing multiple hypotheses.

Finally, it is worth remarking on the connection between these results and
prediction with ``mixtures of experts''
\citep{Arora-Hazan-Kale-multiplicative-weights,prediction-learning-and-games}.
Formally, the role of the negative log-likelihood and of Bayes's rule in this
paper was to provide a loss function and a multiplicative scheme for updating
the weights.  All but one of the main results (Theorem
\ref{thm:generalization-performance}, which bounds Hellinger distance by
Kullback-Leibler divergence) would carry over to multiplicative weight training
using a different loss function, provided the accumulated loss per unit time
converged.

\appendix

\section{Bayesian Updating as Replicator Dynamics}
\label{sec:updating-is-replication}

Replicator dynamics are one of the fundamental models of evolutionary biology;
they represent the effects of natural selection in large populations, without
(in their simplest form) mutation, sex, or other sources of variation.
\citep{Hofbauer-Sigmund-evol-games-and-pop-dyn} provides a thorough discussion.
They also arise as approximations to many other adaptive processes, such as
reinforcement learning
\citep{Borgers-Sarin-reinforcement-and-replicator,Borkar-reinforcement-learning-in-Markovian-games,Sato-JPC-coupled-replicator-equations}.
In this appendix, I show that Bayesian updating also follows the replicator
equation.

We have a set of {\em replicators} --- phenotypes, species, reproductive
strategies, etc. --- indexed by $\parametervalue \in \ParameterSpace$.  The
population density at type $\parametervalue$ is
$\posteriordensity(\parametervalue)$.  We denote by $\phi_t(\parametervalue)$
the {\em fitness} of $\parametervalue$ at time $t$, i.e., the average number of
descendants left by each individual of type $\parametervalue$.  The fitness
function $\phi_t$ may in fact be a function of $\posteriordensity_t$, in which
case it is said to be {\em frequency-dependent}.  Many applications assume the
fitness function to be deterministic, rather than random, and further assume
that it is not an explicit function of $t$, but these restrictions are
inessential.

The discrete-time {\em replicator dynamic}
\citep{Hofbauer-Sigmund-evol-games-and-pop-dyn} is the dynamical system given
by the map
\begin{equation}
\posteriordensity_{t}(\parametervalue)  = \posteriordensity_{t-1}(\parametervalue)\frac{\phi_t(\parametervalue)}{\popavgt{\phi_t}{t}}
\label{eqn:replicator-dynamic}
\end{equation}
where $\popavgt{\phi_t}{t}$ is the population mean fitness at $t$, i.e.,
\[
\popavgt{\phi_t}{t} \equiv \int_{\ParameterSpace}{\phi_t(\parametervalue) d \posteriordensity_t(\parametervalue)}
\]
The effect of these dynamics is to re-weight the population towards replicators
with above-average fitness.

It is immediate that Bayesian updating has the same form as Eq.\
\ref{eqn:replicator-dynamic}, as soon as we identify the distribution of
replicators with the posterior distribution, and the fitness with the
conditional likelihood.  In fact, Bayesian updating is an extra simple case of
the replicator equation, since the fitness function is frequency-independent,
though stochastic.  Updating corresponds to the action of natural selection,
without variation, in a fluctuating environment.  The results in the main text
assume (Assumption \ref{assumption:relative-AEP-holds}) that, despite the
fluctuations, the long-run fitness is nonetheless a determinate function of
$\parametervalue$.  The theorems assert that selection can then be relied upon
to drive the population to the peaks of the long-run fitness function, at the
cost of reducing the diversity of the population, rather as in Fisher's
fundamental theorem of natural selection
\citep{Fisher-genetical,Hofbauer-Sigmund-evol-games-and-pop-dyn}.

\begin{corollary}
  Define the relative fitness $\tilde{\phi}_t(\parametervalue) \equiv
  \CondLike_t(\parametervalue)/\popavgt{\CondLike_t}{t}$.  Under the conditions of
  Theorem \ref{theorem:long-run-fitness}, the time average of the log relative
  fitness converges a.s.
  \begin{equation}
    \frac{1}{t}\sum_{n=1}^{t}{\log{\tilde{\phi}_n(\parametervalue)}} \rightarrow -\Action(\parametervalue) + o(1)
	\label{eqn:cesaro-mean-of-fitness}
  \end{equation}
\end{corollary}

\textsc{Proof:} Unrolling Bayes's rule over multiple observations,
\[
\posteriordensity_t(\parametervalue) = \posteriordensity_0(\parametervalue) \prod_{n=1}^{t}{\tilde{\phi}_n(\parametervalue)}
\]
Take log of both sides, divide through by $t$, and invoke Theorem
\ref{theorem:long-run-fitness}. $\Box$

{\em Remark:} Theorem \ref{theorem:long-run-fitness} implies that
\[
H_t \equiv \left|\log{\posteriordensity_t(\parametervalue)} + t\Action(\parametervalue) \right|
\]
is a.s. $o(t)$.  To strengthen Eq.\ \ref{eqn:cesaro-mean-of-fitness} from
convergence of the time average or Ces{\`a}ro mean to plain convergence
requires forcing $H_t - H_{t-1}$ to be $o(1)$, which it generally isn't.

It is worth noting that Haldane \citep{Haldane-measurement-of-natural-selection}
defined the {\em intensity of selection} on a population as, in the present
notation,
\[
\log{\frac{\posteriordensity_t(\hat{\parametervalue})}{\posteriordensity_0(\hat{\parametervalue})}}
\]
where $\hat{\parametervalue}$ is the ``optimal'' (i.e., most selected-for)
value of $\parametervalue$.  For us, this intensity of selection is just
$\LikelihoodRatio_t(\hat{\parametervalue})/\prioravg{\LikelihoodRatio_t}$ where
$\hat{\parametervalue}$ is the (or a) MLE.

\section{Verification of Assumptions 5--7 for the Example}
\label{appendix:verification}

Since the $X_1^{\infty}$ process is a function of the $S_1^{\infty}$ process,
and the latter is an aperiodic Markov chain, both are $\psi$-mixing (see
\citep{Marton-Shields,Shields-discrete-sample-paths} for the definition of
$\psi$-mixing and demonstrations that aperiodic Markov chains and their
functions are $\psi$-mixing).  Let $\widehat{P^{(k)}_t}$ be the empirical
distribution of sequences of length $k$ obtained from $x_1^t$.  For a Markov
chain of order $k$, the likelihood is a function of $\widehat{P^{(k+1)}_t}$
alone; we will use this and the ergodic properties of the data-generating
process to construct sets on which the time-averaged log-likelihood converges
uniformly.  Doing this will involve constraining both the order of the Markov
chains and their transition probabilities, and gradually relaxing the
constraints.

It will simplify notation if from here on all logarithms are taken to base 2.

Pick $\epsilon > 0$ and let $k(t)$ be an increasing positive-integer-valued
function of $t$, $k(t) \rightarrow \infty$, subject to the limit $k(t) \leq
\frac{\log{t}}{h_{\TrueMeasure}+\epsilon}$, where $h_{\TrueMeasure}$ is the
Shannon entropy rate of $\TrueMeasure$, which direct calculation shows is
$2/3$.  The $\psi$-mixing property of $X_1^{\infty}$ implies \cite[p.\
179]{Shields-discrete-sample-paths} that
\begin{equation}
\TrueMeasure(p_{TV}(\widehat{P^{(k(t)}_t}, \TrueMeasure^{(k(t))}) > \delta)
\leq \frac{\log{t}}{h+\epsilon} 2(n+1)^{t^{\gamma_1}} 2^{-n C_1 \delta^2}
\label{eqn:psi-mixing-bound}
\end{equation}
where $\rho_{TV}$ is total variation distance, $\TrueMeasure^{(k(t))}$ is
$\TrueMeasure$'s restriction to sequences of length $k(t)$, $n = \lfloor t/k(t)
\rfloor -1$, $\gamma_1 = (h_{\TrueMeasure}+\epsilon/2)/(h_{\TrueMeasure} +
\epsilon)$ and $C_1$ is a positive constant specific to $\TrueMeasure$ (the
exact value of which is not important).

The log-likelihood per observation of a Markov chain $\parametervalue \in
\ParameterSpace_k$ is
\[
t^{-1} \log{\modeldensity(x_1^t)} = t^{-1} \log{\modeldensity(x_1^k)} + \sum_{w \in \ObservableSpace^{k}}{\sum_{a \in \ObservableSpace}{\widehat{P^{(k+1)}_t}(wa)\log{\modeldensity(a|w)}}} 
\]
where $\modeldensity(a|w)$ is of course the probability, according to
$\parametervalue$, of producing $a$ after seeing $w$.  By asymptotic
equipartition, this is converging a.s. to its expected value,
$-h_{\TrueMeasure} - \divrate(\parametervalue)$.

Let $z(\parametervalue) = \max_{w,a}{\left|\log{\modeldensity(a|w)}\right|}$.
If $z(\parametervalue) \leq z_0$ and
$\rho_{TV}(\widehat{P^{(k+1)}_t},\TrueMeasure^{(k+1)}) \leq \delta$, then
$t^{-1} \log{\modeldensity(x_1^t)}$ is within $z_0 \delta$ of
$-h_{\TrueMeasure} - \divrate(\parametervalue)$.  Meanwhile, $t^{-1}
\log{\truedensity(x_1^t)}$ is converging a.s. to $-h_{\TrueMeasure}$, and again
\citep{Shields-discrete-sample-paths}
\begin{equation}
\TrueMeasure(|t^{-1} \log{\truedensity(X_1^t)} - h_{\TrueMeasure}| > \delta) \leq q(t,\delta) 2^{-t C_2 \delta}
\label{eqn:exponential-rate-for-entropy-rate}
\end{equation}
for some $C_2 > 0$ and sub-exponential $q(t,\delta)$.  (The details are
unilluminating in the present context and thus skipped.)

Define $G(t,z_0)$ as the set of all Markov models whose order is less than or
equal to $k(t)-1$ and whose log transition probabilities do not exceed $z_0$,
in symbols
\[
G(t,z_0) =\left\{ \parametervalue:z(\parametervalue) \leq z_0\right\} \cap \left(\bigcup_{j=1}^{k(t)-1}{\ParameterSpace_j}\right)
\]
Combining the deviation-probability bounds \ref{eqn:psi-mixing-bound} and
\ref{eqn:exponential-rate-for-entropy-rate}, for all $\parametervalue \in
G(t,z_0)$
\begin{equation}
\TrueMeasure\left(\left|\frac{\log{\LikelihoodRatio_t(\parametervalue)}}{t} - \divrate(\parametervalue)\right| > \delta\right) \leq \frac{\log{t}}{h+\epsilon} 2(n+1)^{t^{\gamma_1}} 2^{-\frac{n C_1 \delta^2}{4z_0}} + q(t,\delta)2^{-\frac{tC_2 \delta}{2}}
\label{eqn:log-likelihood-ratio-deviation-bounds}
\end{equation}
These probabilities are clearly summable as $t\rightarrow\infty$, so by the
Borel-Cantelli lemma, we have uniform almost-sure convergence of $t^{-1}
\log{\LikelihoodRatio_t(\parametervalue)}$ to $-\divrate(\parametervalue)$ for
all $\parametervalue \in G(t,z_0)$.

The sets $G(t,z_0)$ eventually expand to include Markov models of arbitrarily
high order, but maintain a constant bound on the transition probabilities.  To
relax this, let $z_t$ be an increasing function of $t$, $z(t) \rightarrow
\infty$, subject to $z_t \leq C_3 t^{\gamma_2}$ for positive $\gamma_2 <
\gamma_1$.  Then the deviation probabilities remain summable, and for each $t$,
the convergence of $t^{-1} \log{\LikelihoodRatio_t(\parametervalue)}$ is still
uniform on $G(t,z_t)$.  Set $G_t = G(t,z_t)$, and turn to verifying the
remaining assumptions.

Start with Assumption \ref{assumption:exponentially-egorov}; take its items in
reverse order.  So far, the only restriction on the prior $\Posterior_0$ has
been that its support should be the whole of $\ParameterSpace$, and that it
should have the ``Kullback-Leibler rate property'', giving positive weight to
every set $N_{\epsilon} = \left\{\parametervalue: d(\parametervalue) < \epsilon
\right\}$.  This, together with the fact that $\lim_{t}{G_t} =
\ParameterSpace$, means that $\divrate(G_t) \rightarrow
\divrate(\ParameterSpace)$, which is item (3) of the assumption.  The same
argument also delivers Assumption
\ref{assumption:good-sets-are-good-even-for-subsets}.  Item (2), uniform
convergence on each $G_t$, is true by construction.  Finally (for this
assumption), since $\divrate(\ParameterSpace) = 0$, any $\beta > 0$ will do,
and there are certainly probability measures where $\Posterior_0(G_t^c) \leq
\alpha \myexp{-\beta t}$ for some $\alpha, \beta > 0$.  So, Assumption
\ref{assumption:exponentially-egorov} is satisfied.

Only Assumption \ref{assumption:sufficiently-rapid-convergence} remains.  Since
Assumptions
\ref{assumption:measurable-likelihood}--\ref{assumption:relative-AEP-holds}
have already been checked, we can apply Eq.\
\ref{eqn:explicit-upper-bound-on-integrated-likelihood-for-large-t} from the
proof of Lemma \ref{lemma:lim-sup-of-integrated-likelihood-on-good-sets} and
see that, for each fixed $G$ from the sequence of $G_t$, for any $\epsilon >
0$, for all sufficiently large $t$,
\[
t^{-1}\log{\prioravg{G \LikelihoodRatio_t}} \leq -\divrate(G) + \epsilon + t^{-1}\log{\Posterior_0(G)} ~ \mathrm{a.s.}
\]
This shows that $\convtime(G_t, \delta)$ is almost surely finite for all $t$
and $\delta$, but still leaves open the question of whether for every $\delta$
and all sufficiently large $t$, $t \geq \convtime(G_t,\delta)$ (a.s.).
Reformulating a little, the desideratum is that for each $\delta$, with
probability 1, $t < \convtime(G_t, \delta)$ only finitely often.  By the
Borel-Cantelli lemma, this will happen if
$\sum_{t}{\TrueMeasure(\convtime(G_t,\delta) > t)} \leq \infty$.  However, if
$\convtime(G_t,\delta) > t$, it must be equal to some particular $n > t$, so
there is a union bound:
\begin{equation}
  \sum_{t}{\TrueMeasure(\convtime(G_t,\delta) > t)} \leq \sum_{t}{\sum_{n=t+1}^{\infty}{\TrueMeasure\left(\frac{\log{\prioravg{G_t \LikelihoodRatio_n}}}{n} > \delta -\divrate(G_t)\right)}}
\end{equation}
From the proof of Lemma 6 (specifically from Eqs.\
\ref{eqn:integrated-likelihood-on-good-set-broken-up},
\ref{eqn:generating-function-divrate} and
\ref{eqn:generating-function-vs-set-divrate}), we can see that by making $t$
large enough, the only way to have the event $n^{-1}\log{\prioravg{G_t
    \LikelihoodRatio_n}} > \delta -\divrate(G_t)$ is to have
$\left|n^{-1}\log{\LikelihoodRatio_n(\parametervalue)} -
  \divrate(\parametervalue)\right| > \delta/2$ everywhere on a positive-measure
subset of $G_t$.  But we know from Eq.\
\ref{eqn:log-likelihood-ratio-deviation-bounds} not only that the inner sum can
be made arbitrarily small by taking $t$ sufficiently large, but that the whole
double sum is finite.  So $\convtime(G_t,\delta) > t$ only finitely often
(a.s.).

\section*{Acknowledgments}

Portions of this work were done while at the Center for the Study of Complex
Systems at the University of Michigan, supported by a grant from the James
S. McDonnell Foundation, and while at the ``Science et Gastronomie 2003''
workshop organized by the Institut des Syst{\`e}mes Complexes and the
Laboratoire de l'Informatique du Parall{\'e}lisme of the {\'E}cole Normale
Sup{\'e}riere de Lyon.  I thank Michel Morvan and Cristopher Moore for
organizing that workshop; Carl Bergstrom, Anthony Brockwell, Chris Genovese,
Peter Hoff, Giles Hooker, Martin Nilsson Jacobi, Jay Kadane, Kevin T. Kelly,
David Krakauer, Michael Lachmann, Scott E. Page, Mark Schervish, Eric Smith,
Surya Tokdar, and Larry Wasserman for discussion; the anonymous referees for
valuable suggestions; and above all, Kristina Klinkner.

\bibliography{locusts}
\bibliographystyle{acmtrans-ims}

\end{document}